\title{ ~~\\ Value distribution of Ramanujan sums and of cyclotomic polynomial coefficients}
\author{Pieter Moree and Huib Hommerson}
\documentclass[12pt]{article}
\usepackage{amssymb, latexsym, amsfonts}
\def\@ptsize{2}
\newtheorem{Thm}{Theorem}
\newtheorem{Con}{Conjecture}

\newtheorem{Lem}{Lemma}
\newtheorem{cor}{Corollary}

\newtheorem{Prop}{Proposition}
\newcommand{\qed}{\hfill $\Box$}

\begin{document}
\date{}
\maketitle 
{\def\thefootnote{}
}
\begin{abstract}
\noindent The Ramanujan sum $c_n(k)$ and $a_n(k)$, the $k$th coefficient of the $n$th cyclotomic polynomial,
are completely symmetric expressions in terms of primitive $n$th roots of unity. For $k$ fixed we
study the value distribution of $c_n(k)$ (following A. Wintner) and $a_n(k)$ (partly following
H. M\"oller). In particular we disprove a 1970 conjecture of H. 
M\"oller on the average (over $n$) of $a_n(k)$. 
We show that certain symmetric functions in primitive roots considered by the Dence brothers 
are related to the behaviour of $c_{p-1}(k)$ and $a_{p-1}(k)$ as $p$ ranges over the primes and study
their value distribution.\\
\indent This paper is an outgrowth of the M.Sc. thesis project of the second author, carried out under the 
supervision of the first author at the Korteweg-de Vries Institute (University of Amsterdam). Some
of the numerical work described was done by Yves Gallot.
\end{abstract}
\vfil\eject
\tableofcontents
\vfil\eject
\section{Introduction}
We recall from elementary number theory the notion of multiplicative order. If $a$ and $n$ are
coprime natural numbers, then there is a smallest positive integer $k$ such that
$a^k\equiv 1({\rm mod~}n)$, this integer $k$ is the 
{\it multiplicative order} of $a$, ord$_n(a)$, and
divides $\varphi(n)$, where $\varphi$ is Euler's totient (discussed more extensively in the
next section). We say that $a$ is a {\it primitive root} of the multiplicative group 
$(\mathbb Z/n\mathbb Z)^*$ of units of
$\mathbb Z/n\mathbb Z$ if ord$_n(a)=\varphi(n)$. This means that $(\mathbb Z/n\mathbb Z)^*$ is
cyclic and $a$ is a generator. 
Let $p$ be a prime (indeed, throughout we exclusively use the notation $p$ and $q$ for primes).
It is well-known that $(\mathbb Z/p\mathbb Z)^*$ is cyclic and hence it has at least one primitive
root. Indeed, if $g$ is  a primitive root mod $p$, then so are $g^j$ with $j$ and $p-1$ coprime 
and hence there are $\varphi(p-1)$ distinct primitive roots mod $p$. Let us denote them by
$g_1,...,g_t$ with $t=\varphi(p-1)$ and $0<g_i\le p-1$. In their paper \cite{DD} the Dence brothers consider
symmetric functions of the primitive roots of primes. They consider for example
\begin{equation}
\label{eerste}
s_2(p)=\sum_{1\le i<j\le t}g_ig_j
\end{equation}
and show that this quantity, when considered mod $p$, assumes only the values in $\{-1,0,1\}$.
They write in particular (p. 79): `As a matter of distribution, we observe that amongst the
first 100 primes (beginning with $p=5$) the residues $-1,0,1$ of $s_2$ occur in the ratios
$12:59:29$. The question of what these ratios should be in the limit of infinitely many primes
is an interesting one'. The following
table gives the ratios amongst the first $10^j$ primes 
with $2\le j\le 6$ (thus beginning with $p=2$).

\medskip
\centerline{{\bf Table 1:} Value distribution of $s_2(p)$}
\medskip
\begin{center}
\begin{tabular}{|c||c|c|c|c|}
\hline
$\pi(x)$  & $s_2(p)=-1$ & $s_2(p)=0$ & 
$s_2(p)=1$ \\
\hline
$10^2$    & $0.110000$ & $0.610000$ & $0.280000$ \\
\hline
$10^3$ & $0.099000$ & $0.625000$ & $0.276000$ \\
\hline
$10^4$ & $0.093000$ & $0.626100$ & $0.280900$ \\
\hline
$10^5$  & $0.094120$ & $0.627330$ & $0.278550$ \\
\hline
$10^6$  & $0.093939$ & $0.626216$ & $0.279845$ \\

\hline
\end{tabular}
\end{center}
\medskip

\noindent The sum in (\ref{eerste}) is the second elementary totally symmetric function of $g_1,\cdots,g_t$.
For $k\ge 1$ let us consider more generally $s_k(p)$ with $s_k(p)$ the $k$th order totally 
elementary symmetric function in $g_1,\cdots,g_t$.\\

\noindent {\bf Question 1}. Fix $k\ge 1$. As $p$
ranges over the primes, which values are assumed by $s_k(p) ({\rm mod~}p)$ and with
what frequency ?\\
 
\noindent Likewise we consider, for $k\ge 1$, the sum
 $$S_k(p)=\sum_{1\le j\le t}g_j^k.$$ For this
quantity we consider the same question:\\

\noindent {\bf Question 2}. Fix $k\ge 1$. As $p$
ranges over the primes, which values are assumed by $S_k(p) ({\rm mod~}p)$ and with
what frequency ?\\

\noindent One of the motivations of this M.Sc. thesis is to resolve these questions.
Interestingly, we will see that symmetric functions in primitive roots are closely
related with some important objects in number theory; cyclotomic polynomials and Ramanujan sums. 
Both cyclotomic polynomials and Ramanujan sums arise in many contexts and hence we devote 
separate sections to them.\\
\indent We will see that the question of the Dence brothers naturally leads to the study of the
value distribution of Ramanujan sums and of that of
coefficients of cyclotomic polynomials. These two issues have not been well studied in the literature.
Regarding the
value distribution of Ramanujan sums we could only find an almost forgotten paper by the famous 
analyst Aurel Wintner \cite{W}, which we
reconsider in \S \ref{aurelius}. Regarding the value distribution of 
coefficients of cycltomic polynomials, we reconsider a paper by Herbert M\"oller \cite{M} in
\S \ref{vdofcyclotomic} and establish some new results.

\section{Preliminaries}

\subsection{Multiplicative functions}
\label{tweeeen}
An {\it arithmetic function} is a function from the natural numbers to the complex numbers.
A very important subclass of these functions are the so-called {\it multiplicative functions}.
These functions satisfy $f(1)=1$ and
$$f(mn)=f(n)f(m) {\rm ~for~}(m,n)=1.$$
The M\"obius function and Euler totient function are important examples of multiplicative
functions. If $f$ is a multiplicative function, then so is $g(n)$ with
$g(n)=\sum_{d|n}f(d)$. Indeed, more generally, if $f$ and $g$ are two multiplicative functions, so
is $(f\star g)(n)=\sum_{d|n}f(d)g(n/d)$, which is the so called 
{\it Dirichlet convolution} of $f$ and $g$. 
If $n=p_1^{e_1}\cdots p_s^{e_s}$ is the prime factorisation of $n$ and $f$ is
multiplicative, then we have
$$\sum_{d|n}f(d)=\prod_{j=1}^s\left(1+f(p)+\cdots+f(p^{e_j})\right).$$
A closely related formula is the following one; if $\sum_{d=1}^{\infty}f(d)$
is absolutely convergent and $f$ is multiplicative, then
\begin{equation}
\label{euleridentity}
\sum_{d=1}^{\infty}f(d)=\prod_p\left(\sum_{j=0}^{\infty}f(p^j)\right).
\end{equation}
The latter identity is known as {\it Euler's product identity}.\\
\indent If $f$ is multiplicative, then it is an easy observation that
\begin{equation}
\label{spreiding}
f(m)f(n)=f((m,n))f([m,n])
\end{equation}
for all integers $m$ and $n$, where $(m,n)$ denotes the greatest common 
divisor of $m$ and $n$ and $[m,n]$ the lowest common multiple. Likewise
we define $(m_1,m_2,\dots,m_s)$ and $[m_1,m_2,\dots,m_s]$, 
were we put $(m_1)=1$ and $[m_1]=m_1$.
Conversely,
a function that satisfies the functional equation (\ref{spreiding}) is
said to be {\it semi-multiplicative}. It can be shown \cite{R} that if $f$
is semi-multplicative, there exists a non-zero constant $c$, a positive integer
$a$, and a multiplicative function  $f_1$ such that
$$f(n)=\cases{cf_1(n/a) &if $a\nmid n$;\cr
0 & otherwise.}$$
It is easy to see that if $f$ is semi-multiplicative and $k$ is a constant,
then also $f(k/(k,n))$ is semi-multiplicative in $n$.

\subsection{The M\"obius function}
Let $n$ be an integer having prime divisors $p_1,\cdots,p_s$.
Then the M\"obius $\mu$ function is defined as follows:
$$\mu(n)=\cases{0 &if $p^2|n$ for some prime $p$;\cr 
(-1)^s &otherwise.}$$
Note that $\mu(1)=1$ and that $\mu(n)^2=1$ iff $n$ is squarefree. 
The M\"obius function satisfies the following identity:
\begin{equation}
\label{een}
\sum_{d|n}\mu(d)=\cases{1 &if $n=1$;\cr
0 &otherwise.}
\end{equation}
\noindent An integer $n$ is said to be $k${\it th powerfree} if $n$ is not divisible by
a $k$th power of some integer $>1$. Later we will need the following identity.
\begin{Prop}
\label{kmoebius}
For any natural number $k$ we have
$$\sum_{d^k|n}\mu(d)=\cases{1 &if $n$ is $k$th power-free;\cr
0 &otherwise.}$$
\end{Prop}
First proof. Using the multiplicativity of the M\"obius-function one
sees that the left hand side is a multiplicative function of $n$. It
thus suffices to evaluate it when $n$ is a prime power. Now
$$\sum_{d^k|p^{\alpha}}\mu(d)=\cases{1 &if $\alpha\le k-1$;\cr
0 &otherwise,}$$
from which the result follows. \qed\\
Second proof. Write $n=m^kr$, with the number $r$ $k$th power free. Then the sum under consideration
equals $\sum_{d|m}\mu(d)$. The result then follows by (\ref{een}). \qed \\

\noindent An important formula involving the M\"obius function is the celebrated {\it M\"obius inversion
formula}. It states that $f(n)=\sum_{d|n}g(d)$ iff $g(n)=\sum_{d|n}\mu(d)f(n/d)$. The M\"obius function
often arises in combinatorial problems, we will see several examples in this M.Sc. thesis. In Section
\ref{zes} some auxiliary functions involving the M\"obius function are considered.

\subsection{The Euler totient function}

The M\"obius function often arises in
combinatorial problems where so-called inclusion and exclusion is being
used in counting. We demonstrate this by deriving an identity for
Euler's totient function $\varphi(n)$, which is defined as 
$$\varphi(n)=\sum_{1\le j\le n,~(j,n)=1}1.$$
\begin{Prop}
\label{eulerphi}
Let $N(d)$ denote the number of integers $1\le m\le n$ that are divisible by $d$.
Then
$$\varphi(n)=\sum_{d|n}\mu(d)N(d)=n\sum_{d|n}{\mu(d)\over d}=n\prod_{p|n}\left(1-{1\over p}\right).$$
\end{Prop}
{\it Proof}. Suppose $1\le m\le n$ is an integer. In the expression
$\sum_{d|n}\mu(d)N(d)$, the integer $m$ 
is counted in those $N(d)$ for which both $d|m$ and $d|n$.
It is counted with weight $\sum_{d|n,~d|m}\mu(d).$
But, using  (\ref{een}), we see that
$$\sum_{d|n,~d|m}\mu(d)=\sum_{d|(n,m)}\mu(d)=\cases{1 &if $(n,m)=1$;\cr
0 &otherwise}.$$
Thus $\phi(n)=\sum_{d|n}\mu(d)N(d)$. In case $d|n$, then $N(d)=n/d$.
On using that $\mu$ is a multiplicative function, the result then follows. \qed\\

\noindent From Proposition \ref{eulerphi} we easily infer that $\phi$ is a multiplicative 
function. Alternatively this can be easily deduced on using the Chinese remainder theorem.\\ 
\indent For later use we note that by Proposition \ref{eulerphi} we have
\begin{equation}
\label{varquotient}
{\varphi(\delta n)\over \varphi(n)}=\delta{\prod_{p|n\delta}(1-{1\over p})\over \prod_{p|n}(1-{1\over p})}
=\delta\prod_{p|\delta\atop p\nmid n}(1-{1\over p})=\delta{\varphi(\delta_1)\over \delta_1},
\end{equation}
where $\delta_1=\prod_{p|\delta,~p\nmid n}p$. Note that as $n$ runs over
all integers, $\delta_1$ will run over all squarefree divisors of $\delta$.
In particular, (\ref{varquotient}) implies that
$$
{\varphi(2n)\over \varphi(n)}=
\cases{
2 &if $n$ is even;\cr
1 &otherwise.}
$$

\subsection{Roots of unity}
An $n$th root of unity is  a solution of $z^n=1$ in $\mathbb C$. Note
there are precisely $n$ solutions; $e^{2\pi i/n},\cdots,e^{2\pi in/n}$, with,
of course, $e^{2\pi in/n}=1$.
Instead of $e^{2\pi ik/n}$ one often writes $\zeta_n^k$ with
$\zeta_n=e^{2\pi i/n}$. Indeed, throughout this M.Sc. thesis we adopt the
notation
$$\zeta_n=e^{2\pi i\over n}.$$
An $n$th root of unity is said to be a {\it primitive} root of unity if
it is of the form $\zeta_n^k$ with $k$ and $n$ coprime (which we will
denote by $(k,n)=1$). Such a primitive root of unity has the property that it
does not satisfy an identity of the form $z^m=1$ with $m$ a divisor $<n$ of
$n$.\\

\noindent {\tt Example}. For $n=6$, $\zeta_6$ and $\zeta_6^5={\bar \zeta_6}$ are
the primitive roots of unity. Note that $\zeta_6^2$ is not a primitive
$6$th order root of unity as it satisfies the identity $z^3=1$.\\ 

\noindent Note that there are $\varphi(n)$ primitive $n$th roots of unity. The
roots of unity form a cyclic group of order $n$. The primitive roots of unity
correspond with the generators of this group.

\section{Ramanujan sums}
The Ramanujan sum $c_n(m)$ is defined by
$$c_n(m)=\sum_{1\le k\le n\atop (k,n)=1}e^{2\pi i mk\over n}=\sum_{1\le k\le n\atop (k,n)=1}\zeta_n^{mk}.$$
Although Ramanujan was not the first to work with Ramanujan sums, he was the first tor realize their importance
and use them consistently (especially in the theory of represenation of numbers as sum of squares, see 
e.g \cite[Chapter IX]{Hardy}).
The next proposition lists some of the important properties of Ramanujan sums. 
(These properties are all known or well-known except perhaps the formula in part 6, for which we do not
have a reference.)
By $\nu_p(n)$ we will
denote the exponent of $p$ in $n$, that is we have $\nu_p(n)=r$ iff $p^r|n$ and $p^{r+1}\nmid n$.
\begin{Prop}
\label{basicramanujan}
\item{\rm 1)} We have
$$c_n(m)=\sum_{d|(n,m)}d\mu({n\over d}).$$
\item{\rm 2)} We have
$$c_n(m)=\mu\left({n\over (n,m)}\right){\varphi(n)\over \varphi({n\over (n,m)})}.$$
\item{\rm 3)} We have $-(n,m)\le c_n(m)\le (n,m)$ and $-\varphi(n)\le c_n(m)\le \varphi(n)$.
\item{\rm 4)} We have $c_n(m)\in \mathbb Z$.
\item{\rm 5)} We have $c_{n_1n_2}(m)=c_{n_1}(m)c_{n_2}(m)$ if $(n_1,n_2)=1$; i.e. $c_n(m)$ is
multiplicative with respect to $n$.
\item{\rm 6)} The function $c_n(m)$ is 
semi-multiplicative in $m$. For fixed squarefree $n$, the function $\mu(n)c_n(m)$ is multiplicative in $m$. For
arbitrary natural numbers $n$ and $m$ we have
$$c_n(m)=\prod_{p|m}\mu\left({p^{\nu_p(n)}\over (n,p^{\nu_p(m)})}\right){\varphi(n)\over \varphi({n\over (n,p^{\nu_p(m)})})}.$$
\item{\rm 7)} The following orthogonality relations hold (when $r_1|r$ and $r_2|r$):
$${1\over r}\sum_{m=1}^r c_{r_1}(m)c_{r_2}(m)=\cases{0 &if $r_1\ne r_2$;\cr
\varphi(r_1) &if $r_1=r_2$.}$$
\end{Prop}
{\it Proof}. 1) Consider $g(n):=\sum_{1\le k\le n}e^{2\pi i mk\over n}$. Clearly this expression equals $n$ if $n|k$ and
zero otherwise. Expressing $g(n)$ in terms of Ramanujan sums we obtain
$g(n)=\sum_{d|n}c_{n/d}(m)$, which by M\"obius inversion gives 
$c_n(m)=\sum_{d|n}g(d)\mu(n/d)$. On using that $g(d)=d$ if $d|m$ and vanishes otherwise, the stated
formula follows.\\
2) This result, due to O. H\"older, follows from property 1 on using (\ref{varquotient}) (see e.g. 
Exercise 1.1.14 of \cite{Murty}).\\ 
3) This is a consequence of property 2 and formula (\ref{varquotient}).\\
4) First proof. Immediate by property 1. Second proof. The Ramanujan sum can be considered as an element of the field extension
$\mathbb Q(\zeta_n)$. It is invariant under each of the automorphisms $\zeta_n\rightarrow \zeta_n^k$ with
$k$ coprime to $n$. Hence it is in the fixed field of $\mathbb Q(\zeta_n)$, which is
$\mathbb Q$, since $\mathbb Q(\zeta_n):\mathbb Q$ is Galois. Since $\zeta_n$ is an algebraic integer, the
Ramanujan sum must even be an integer.\\
5) First proof. We have
$$c_{n_1}(m)c_{n_2}(m)=\sum_{1\le k_1\le n_1\atop (k_1,n_1)=1}e^{2\pi i mk_1\over n_1}
\sum_{1\le k_2\le n_2\atop (k_2,n_2)=1}e^{2\pi i mk_2\over n_2}
=\sum_{k_1,k_2} e^{2\pi i m(k_1n_2+k_2n_1)\over n_1n_2}=c_{n_1n_2}(m),$$
where we use the observation that the set of congruences classes of the form $k_1n_2+k_2n_1$, with 
$1\le k_j\le n_j$ and $(k_j,n_j)=1$ for $j=1$ and $j=2$, consists of $\varphi(n_1n_2)$ distinct congruence classes mod $n_1n_2$, all of them
coprime with $n_1n_2$.\\
Second proof. The assumptions on $n_1$ and $n_2$ imply that $(n_1n_2,m)=(n_1,m)(n_2,m)$. Now note that
$$c_{n_1n_2}(m)=\sum_{d|(n_1n_2,m)}d\mu({n_1n_2\over d})=\sum_{d_1|(n_1,m)\atop d_2|(n_2,m)}d_1d_2\mu({n_1\over d_1}
{n_2\over d_2})=c_{n_1}(m)c_{n_2}(m),$$
where in the last step we use the multiplicativity of $\mu$.\\
6) The semi-multiplicativity is immediate from what has been said in Section \ref{tweeeen} and property 2. 
The second assertion is an easy consequence of property 2 and the observation that if $f$ is 
multiplicative, then $f(n)/f(n/(n,k))$ is a multiplicative function in $k$. The proof of the
third assertio uses in addition to the previous proof ingredients the multiplicativity of $\mu$.\\
7) For a proof see e.g. p. 17 of \cite{SS}. \qed\\

\noindent An arithmetic function $f$ is said to be an ${\it even}$ function of
$(n,r)$ if $f((n,r),r)=f(n,r)$ for all $n$ (for a survey see \cite{C}). By Property 2 
$c_n(m)=c_n((n,m))$ and hence Ramanujan sums are even. By Property 4 we can alternatively
define $c_n(m)=\sum_{1\le k\le n,~(k,n)=1}\cos(2\pi mk/n)$.\\
\indent An optimality property of Ramanujan sums was discovered by Bachman \cite{B0}. If $r\ge 1$ is
any real number, he showed that for any sequence of real numbers $b_k$ we have
$$\sum_{m=1}^n|\sum_{1\le k\le n\atop (k,n)=1}b_ke^{2\pi i m k/n}|^r \ge
\left(|\sum_{1\le k\le n,~(k,n)=1}b_k|\over \varphi(n)\right)^r \sum_{m=1}^{n}|c_n(m)|^r.$$
It follows that if we consider the infimum of the left handside over all sequences $b_k$
with $b_k$ this is assumed in case $b_k=1$ for all $k$.\\ 
\indent For some further properties of Ramanujan
sums, such as for example the Brauer-Rademacher identity, we refer to Chapter 2 of \cite{Mc}.

\subsection{Density and order of growth}
If $S$ is a set of natural numbers, by $\delta(S)$ we denote the limit of $x^{-1}\sum_{n\le x,~n\in S}1$ 
(as $x$ tends to infinity) if this
exists. Let $\pi(x)$ denote the number of primes $p\le x$. 
Recall that the prime number theorem asserts that asymptotically $\pi(x)\sim x/\log x$ (for a much stronger
version see Lemma \ref{siegelwalfisz} below).
If $S$ is a set of 
prime numbers, by $\delta(S)$ we denote the limit of $\pi(x)^{-1}\sum_{p\le x,~p\in S}1$ (as $x$ tends
to infinity) if this
exists.\\
\indent If $f$ and $g$ are two functions defined on a set $S$, then we shall write $f(x)=O(g(x))$ (we 
say $f$ is of {\it order} $g$) if
the ratio $|f(x)/g(x)|$ is bounded for all $x\in S$. If the ratio $f(x)/g(x)$ tends to zero for $x$
tending to a specified value $x_0$ (which may be infinite), then we shall write $f(x)=o(g(x))$.
The symbols $O(\cdot)$ and $o(\cdot)$ are commonly called {\it Landau symbols}.

\subsection{Ramanujan expansions of arithmetic functions}
Ramaujan sums play a key role in the theory of {\it Ramanujan expansions}. For completeness
we discuss this topic, though this will be not needed in the rest of this thesis.\\
\indent Let $$M(f)=\lim_{x\rightarrow \infty}{1\over x}\sum_{n\le x}f(n).$$
denote the mean-value of $f:\mathbb N\rightarrow \mathbb C$ (if it exists). The Ramanujan
sums have the following orthogonality property which is easily deduced from Property 7 of Proposition 
\ref{basicramanujan}:
$$M(c_r\cdot c_s)=\cases{\varphi(r) &if $r=s$;\cr 0 &otherwise.}$$
On certain spaces of arithmetic functions there is an inner-product $\langle f\cdot g\rangle=M(f\cdot {\bar g})$
and this product suggests {\it Ramanujan expansions} $f\sim \sum_{r\ge 1}a_r(f)c_r$ for
arithmetic functions $f$ in a suitable Hilbert space with coefficients $a_r(f)=\langle f\cdot c_r\rangle /\varphi(r)$.
Indeed, Ramanujan showed that for example
$${\varphi(n)\over n}={6\over \pi^2}\sum_{r\ge 1}{\mu(r)c_r(n)\over \prod_{p|r}(p^2-1)},$$
$$\sigma_s(n)=\sum_{d|n}d^s={n^s\over \zeta(s)}\sum_{r=1}^{\infty}{\mu(r)c_r(n)\over k^s\prod_{p|k}(1-1/p^s)}~~
({\rm Re}(s)>1),$$
and
$$r(n)=\pi \sum_{r\ge 1}{(-1)^{r-1}c_{2r-1}(n)\over 2r-1},$$
where $r(n)$ denotes the number of representations of $n$ as a sum of two squares and $\zeta(s)$ 
as usual the Riemann zeta function (for
Re$(s)>1$ we have $\zeta(s)=\sum_{n=1}^{\infty}n^{-s}$ and for all other $s\ne 1$ in the complex plane
$\zeta(s)$ can be uniquely
defined by analytic continuation). Rearick \cite{R} proved that if $F$ is an arithmetical function and
an associated function $f$ is defined by $F(n)=\sum_{d|n}df(d)$ such that $\sum_k \sum_n |f(kn)|<\infty$, then
we have, with $a(n)=\sum_k f(kn)$, that the series $\sum_k a(k)c_k(n)$ converges absolutely for each $n$ and we have
$F(n)=\sum_k a(k)c_k(n)$.
For a readable survey on Ramanujan expansions, see e.g. \cite{Schwarz}, for more detailed information
see \cite{SS}.

\subsection{Ramanujan sums and Dirichlet series}
\label{ramadirichlet}
Given an arithmetic function $f$ an important quantity in analytic number
theory is the associated generating series, $\sum_{n=1}^{\infty}f(n)/n^s$, which is usually named after
Dirichlet. If we take $f(n)=c_n(m)$, we obtain
$$R_m(s)=\sum_{n=1}^{\infty}{c_n(m)\over n^s}.$$
By property 1 of Proposition \ref{basicramanujan} we have $c_n(m)=\sum_{dr=n,~d|m}\mu(r)d$ and hence
$${c_n(m)\over n^s}=\sum_{dr=n,~d|m}{\mu(r)\over r^{s}}d^{1-s}$$
and so
$$R_m(s)=\sum_{n=1}^{\infty}\sum_{dr=n,~d|m}{\mu(r)\over r^{s}}d^{1-s}=
\sum_{d|m}d^{1-s}\sum_{r=1}^{\infty}{\mu(r)\over r^s}
={\sigma_{1-s}(m)\over \zeta(s)}={m^{1-s}\sigma_{s-1}(m)\over \zeta(s)}.$$
(This proof is due to Estermann, for Ramanujan's (longer) proof see \cite[p. 140]{Hardy}).
Taking $s=2$ we obtain in particular
$$\sigma_1(m)=\sum_{d|m}d={\pi^2\over 6}m\left(1+{(-1)^m\over 2^2}
+{2\cos(2m\pi/3)\over 3^2}+{2\cos(m\pi/4)\over 4^2}+\cdots \right),$$
which shows in a striking way the oscillation of $\sigma_1(m)$ about its `average' $\pi^2 m/6$.
(It is not difficult to show that $\sum_{m\le x}\sigma_1(m)\sim {\pi^2\over 6}\sum_{m\le x}m\sim {\pi^2\over 12}x^2$
as $x$ tends to infinity:
$$\sum_{m\le x}\sigma_1(m)=\sum_{rd\le x}d={1\over 2}\sum_{r\le x}[{x\over r}]
([{x\over r}]+1)={1\over 2}\sum_{r\le x}{x^2\over r^2}+O\left(x\sum_{r\le x}{1\over r}\right),$$
for which it follows that $\sum_{m\le x}\sigma_1(m)={\pi^2x^2/12}+O(x\log x)$.)\\
\indent Similarly one can consider $R_n(s)=\sum_{m=1}^{\infty}c_n(m)/m^s$. Invoking property 1 of 
Proposition \ref{basicramanujan} again, one
easily obtains that $R_n(s)=\zeta(s)\sum_{d|n}\mu(n/d)d^{1-s}$.\\
\indent In \S \ref{aurelius} the Dirichlet series $f_m^{(j)}(s)=\sum_{n=1}^{\infty}c_n(m)^j/n^s$ features (with
$j$ a non-negative integer).
Define
\begin{equation}
\label{eerstewinter}
F_m^{(j)}(s)=\prod_{p|m}{\left(1+\sum_{k=1}^{\nu_p(m)}{(p^k-p^{k-1})^j\over p^{ks}}+{(-p^{\nu_p(m)})^j\over
p^{(\nu_p(m)+1)s}}\right) \over 1+{(-1)^j\over p^s}}.
\end{equation}
Using Properties 5 and 2 of Proposition \ref{basicramanujan} and Euler's product identity, one infers that, 
for Re$(s)>1$, 
$$f_m^{(j)}(s)=F_m^{(j)}(s)\prod_p(1+{(-1)^j\over p^s}).$$
It follows that
\begin{equation}
\label{wrelatie}
f_m^{2j}(s)=F_m^{(2j)}(s)\zeta(s)/\zeta(2s){\rm ~and~}
f_m^{2j-1}(s)=F_m^{(2j-1)}(s)/\zeta(s).
\end{equation}

\section{Cyclotomic polynomials}
The $n$th cyclotomic polynomial, $\Phi_n(x)$, is defined by
\begin{equation}
\label{definitie}
\Phi_n(X)=\prod_{1\le j\le n\atop (j,n)=1}(X-\zeta_n^j).
\end{equation}
Note that the roots of $\Phi_n(X)$ are precisely the primitive $n$th roots
of unity. Clearly $\Phi_n(X)$ is a monic polynomial of degree $\varphi(n)$.
We write $$\Phi_n(X)=\sum_{k=0}^{\phi(n)}a_n(k)x^k.$$ It is the coefficients
$a_n(k)$ that later will have our special attention. Using induction with
respect to $n$ it is not difficult to show that $a_n(k)\in \mathbb Z$ 
(see e.g. \cite[Exercise 1.5.24]{Murty}).
Note that
$$X^n-1=\prod_{d|n}\prod_{1\le j\le n\atop (j,n)=d}(X-\zeta_n^j)=\prod_{d|n}\Phi_{n\over d}(X)=
\prod_{d|n}\Phi_d(X).$$
By applying M\"obius inversion one infers from this that
\begin{equation}
\label{basiccyclo}
\Phi_n(X)=\prod_{d|n}(X^d-1)^{\mu({n\over d})}.
\end{equation}
Thus, for example,
$$\Phi_p(X)={X^p-1\over X-1}=X^{p-1}+X^{p-2}+\cdots+X+1.$$
Furthermore, it is known that $\Phi_{2n}(X)=\Phi_n(-X)$ if $n>1$ is odd and
$\Phi_{n}(X)=\Phi_{\gamma(n)}(X^{n/\gamma(n)}),$ where $\gamma(n)=\prod_{p|n}p$ is the
{\it squarefree kernel} of $n$. The latter two properties are the reason that some authors in
this area restrict themselves to $n$ squarefree and odd. A further important property is
that for $n>1$ we have $X^{\varphi(n)}\Phi_n(1/X)=\Phi_n(X)$. 
It is not difficult to infer these three
properties from (\ref{basiccyclo}).
In terms of the coefficients
these properties imply (respectively):
\begin{equation}
\label{naarkwadraatvrijekern}
a_n(k)=\cases{a_{\gamma(n)}({k\gamma(n)\over n}) &if ${n\over \gamma(n)}|k$;\cr
0 & otherwise,}
\end{equation}
\begin{equation}
\label{verdubbeling}
a_{2n}(k)=(-1)^ka_n(k)~{\rm ~for~}n>1,~2\nmid n;
\end{equation}
$$a_n(k)=a_n(\varphi(n)-k)~{\rm ~for~}n>1.$$ 
It can be shown that $\Phi_n(X)$ is irreducible over $\mathbb Q$. This has as
consequence that the degree of the cyclotomic field $\mathbb Q(\zeta_n)$ over $\mathbb Q$
equals $\varphi(n)$. Over finite fields the cyclotomic fields factor in general. Indeed,
over $\mathbb F_q$, $\Phi_d(X)$ has $\phi(d)/{\rm ord}_q(d)$ distinct irreducible factors if
$(q,n)=1$. If $p\nmid m$, then $p|\Phi_m(a)$ if and only if the order of $a({\rm mod~}p)$ 
is $m$ (Proposition \ref{flauwflauw} below). This can be used to infer the infinitude of primes $p\equiv 1({\rm mod~}m)$ 
\cite[Exercise 1.5.30]{Murty}.\\
\indent Note that the coefficient, $a_n(\varphi(n)-1)$, of $X^{\varphi(n)-1}$ equals
$$-\sum_{1\le j\le n\atop (j,n)=1}\zeta_n^j=-c_n(1)=-\mu(n),$$
where the last equality follows from Property 5 of Proposition \ref{basicramanujan}.
Since $a_n(\varphi(n)-1)=a_n(1)$ for $n>1$, we conclude that $a_n(1)=-\mu(n)$.
Another connection between cyclotomic polynomials and Ramanujan sums is given
by the following result due to Nicol \cite{Nicol}.
\begin{Prop}
We have
$$\Phi_n(X)=\exp\left(-\sum_{m=1}^{\infty}{c_n(m)\over m}X^m\right)$$
and
$$\sum_{m=1}^n c_n(m)X^{m-1}=(X^n-1){\Phi'_n(X)\over \Phi_n(X)}.$$
\end{Prop}
Using cyclotomic polynomials we can give an easy proof of Property 1 of Proposition \ref{basicramanujan}.
We take the logarithmic of the right hand side of (\ref{definitie}) and
compare its Taylor series with that obtained on taking the logarithm of the right hand side
of (\ref{basiccyclo}). The details are left to the reader.

\subsection{On the size of the coefficients of cyclotomic polynomials}

The size of the coefficients of cyclotomic polynomials is a much researched topic. This presumably
stems from wonder over the fact that so often the coefficients are in $\{0,\pm 1 \}$. Indeed,
only for $n\ge 105$ some coefficients outside this range appear. For example, we 
have $a_{105}(7)=-2$. 
The amazement over the smallness of $a_n(m)$ was eloquently worded by D. Lehmer \cite{DL}
who wrote: `The smallness of $a_n(m)$ would appear to be one of the fundamental conspiracies
of the primitive $n$th roots of unity. When one considers that $a_n(m)$ is a sum of
$({\phi(n)\atop m})$ unit vectors (for example $73629072$ in the case of $n=105$ and
$m=7$) one realizes the extent of the cancellation that takes place'.\\ 
\indent Migotti showed in 1883 that $a_{pq}(i) \in \{0,\pm 1\}$, with
$p$ and $q$ odd primes. On the other hand Emma Lehmer showed that $a_{pqr}(i)$ can be
arbitrarily large (with $p,q$ and $r$ odd primes).
Sister Marion Beiter put forward in 1971 the conjecture that in 
fact $a_{pqr}(i)\le (p+1)/2$ for all $i$ and for any $p<q<r$ and
that this bound is best possible. This conjecture remains unsolved, for some 
recent progress see \cite{B1}.\\
\indent Let us put $B(k)=\max_{n\ge 1}|a_n(k)|$. 
For small values of $k$, M\"oller \cite{HM} gave a method to compute $B(k)$. 
The next table (taken from \cite{HM}) gives some values of $B(k)$.\\
\vfil\eject

\centerline{{\bf Table 2:} Value of $B(k)$ for $1\le k\le 30$}
\begin{center}
\begin{tabular}{|c|c|c|c|c|c|c|c|c|c|c|c|c|c|c|c|c|}
\hline
$k$  & $1$ & $2$ & $3$ & $4$ & $5$ & $6$ & $7$ & $8$ & $9$ & $10$ & $11$ & $12$ & $13$ & $14$ & $15$ \\
\hline
$B(k)$  & $1$ & ${1}$ & ${1}$ & ${1}$ & ${1}$ & 
$1$ & $2$ & $1$ & ${1}$ & $1$ & $2$ & $1$ & $2$ & $2$ & $2$\\
\hline
$k$ & $16$ & $17$ & $18$ & $19$ & $20$ & $21$ & $22$ & $23$ & $24$ & $25$ & $26$ & $27$ & $28$ & $29$ & $30$\\
\hline
$B(k)$ & $2$ & $3$ & $3$ & $3$ & $3$ & $3$ &  $3$ & $4$ & $3$ & $3$ & $3$ & $3$ & $4$ & $4$ & $5$ \\
\hline
\end{tabular}
\end{center}
\medskip
The table suggests that $B(k)\ge k$ for every $k\ge 1$, this is however very far from
the case: Bachman \cite{B} extended work of several earlier authors and showed that
\begin{equation}
\label{gennady}
\log B(k)=C_0{\sqrt{k}\over (\log k)^{1/4}}\left(1+O\left({\log \log k\over \sqrt{\log k}}\right)\right).
\end{equation}
In the opposite direction, let us put $A(n)=\max_{m}|a_n(m)|$. Erd\"os has shown that there
exists a $c>0$ and that there are infinitely many $n$ such that
$$\log A(n) \gg \exp\left({c\log n\over \log \log n}\right).$$ On the other hand it is known
that
$$\log A(n) < \exp\left((\log 2+o(1)){\log n\over \log \log n}\right),$$
where the constant $\log 2$ is best possible. So the conspiration of the primitive roots of
unity is not always that efficient !\\
\indent Jiro Suzuki \cite{Suzuki} proved the following beautiful result, which we present
with a slightly shortened proof. Proposition \ref{suzukieven}, using a slightly more
involved argument, 
strengthens Suzuki's result.
\begin{Prop}
\label{z}
We have $\{a_{n}(k):n,k\in \mathbb N\}=\mathbb Z$.
\end{Prop}
{\it Proof}.  
Given any integer $s\ge 2$ it is a consequence of the prime number theorem (cf. the
proof of Proposition \ref{bertie}) that
it is possible to find primes $2<p_1<p_2<\cdots <p_s$ such that $p_1+p_2>p_s$. 
Let $q$ be any prime $>p_s$. If $s$ is even, let $m=p_1\cdots p_sq$, otherwise
let $m=p_1\cdots p_s$. We claim that
\begin{equation}
\label{suziall}
\Phi_m(X)=\sum_{j=0}^{p_s-1}a_m(j)X^j+(1-s)X^{p_s} ~({\rm mod~}X^{p_s+1}).
\end{equation}
The claim shows that $a_m(p_s)=1-s$. By (\ref{verdubbeling}) we have $a_{2m}(p_s)=s-1$. On
noting that $a_4(1)=0$ the result follows. We next prove (\ref{suziall}). 
By (\ref{thanga}) and
the observation that $\{d:d|m,~d<p_s+1\}=\{1,p_1,\dots,p_s\}$ and
that for $i\ne j$, $p_i+p_j\ge p_1+p_2>p_s$, we infer 
 that mod $X^{p_s+1}$ we
have
\begin{eqnarray}
\Phi_m(X)&\equiv &{(1-X^{p_1})\cdots (1-X^{p_s})/(1-X)}~({\rm mod~}X^{p_s+1})\cr
&\equiv & (1+X+\cdots+X^{p_s-1})(1-X^{p_1}-\cdots-X^{p_{s-1}})~({\rm mod~}X^{p_s+1})\cr
&\equiv &\sum_{j=0}^{p_s-1}a_m(j)X^j+(1-s)X^{p_s} ~({\rm mod~}X^{p_s+1}).\nonumber
\end{eqnarray}
This concludes the proof. \qed\\

\noindent The basic idea of this proof seems to have been first formulated by I. Schur in a letter
to E. Landau. He used it to show that $a_n(k)$ can be arbitrarily large.\\ 
\indent A very readable survey on properties of coefficients of cyclotomic polynomials is 
provided by \cite{Thanga}.\\ 
\indent In this thesis we are especially interested in how often certain values are assumed by
$a_n(k)$ for small fixed $k$ (\S \ref{vdofcyclotomic}), where $n$ either runs over all integers $n$ or
over the numbers of the form $p-1$ with $p$ a prime (\S \ref{vdofsk}).

\section{Primitive roots}
\label{primitiverOOts}
We already stated that for every prime $p$ there exists a primitive root $g$ mod $p$ and
that all other primitive roots mod $p$ are of the form $g^j$ with $j$ coprime to $p-1$.
Since $g,g^2,...,g^{p-1}$ are all distinct mod $p$, there are precisely $\varphi(p-1)$ distinct
primitive roots mod $p$.\\
\indent Emil Artin conjectured in 1927 that for every integer $g$ not equal to $-1$ or a square, there
are infinitely many primes $p$ such that $g$ is  a primitive root. More precisely, he conjectured
that the set of primes $p$ for which $g$ is a primitive root has a density which should equal
a rational multiple of $A$, the Artin constant, which is defined as
$$A=\prod_p(1-{1\over p(p-1)})=0.3739558136\cdots$$
Artin's conjecture was proved in 1967 under the assumption of the Generalized Riemann Hypothesis (GRH)
by C. Hooley. In our consideration of the problem stated in the introduction, the constant $A$ will
also arise. Artin's constant (and indeed many similar ones) can be evaluated with many decimals of precision
by expanding it in terms of values of the Riemann zeta-function evaluated at integers exceeding one, 
see e.g. \cite{Mor}. We return to the Artin primitive root conjecture in Section \ref{return}.

\section{Some sums involving the M\"obius function}
\label{zes}
Two important auxiliary quantities that will arise are, for given positive integers $k,r$ and $j$, are the
sums
\begin{equation}
\label{twosums}
\sum_{n\le x\atop (n,r)=1}\mu(n)^j {\rm ~and~}\sum_{p\le x\atop (p-1,k)=r}\mu(p-1)^j.
\end{equation}
A number of the form $p-t$ with $t$ fixed is said to be a {\it shifted prime}.
Note that we can restrict to the case where $j=1$ and $j=2$.
The case $j=2$ can be dealt with by fairly elementary analytic number theory. The case
$j=1$ is rather harder, in the case of shifted primes providing a non-trivial estimate
would be considered a major achievement by the experts we consulted.

\subsection{Counting squarefree integers}
\noindent In this section we are concerned with estimating the first sum appearing in (\ref{twosums}).
\begin{Lem}
\label{mobiuseen}
Let $r\ge 1$ be an integer. We have
$$\sum_{m\le x\atop (m,r)=1}\mu(m)^2={6x\over \pi^2\prod_{p|r}(1+{1\over p})}+O(\sqrt{x}\varphi(r)),$$
where the implied constant is absolute.
\end{Lem}
{\it Proof}. We have, by inclusion and exclusion,
$$\sum_{m\le x\atop (m,r)=1}\mu(m)^2=\sum_{d\le \sqrt{x}\atop (d,r)=1}\mu(d)A_r({x\over d^2}),$$
where $A_r(x)$ denotes the number of integers $n\le x$ that are coprime with $r$. Note that
$$[{x\over r}]\varphi(r)\le A_r(x)\le [{x\over r}]\varphi(r)+\varphi(r)$$
and hence $A_r(x)=\varphi(r)x/r+O(\varphi(r))$. On using the latter estimate we obtain
\begin{eqnarray}
\sum_{m\le x\atop (m,r)=1}\mu(m)^2&=&x{\varphi(r)\over r}\sum_{d\le \sqrt{x}\atop (d,r)=1}{\mu(d)\over d^2}
+O(\sqrt{x}\varphi(r)).\nonumber\\
&=&x{\varphi(r)\over r}\sum_{(d,r)=1}^{\infty}{\mu(d)\over d^2}+O(\sqrt{x}\varphi(r)).\nonumber\\
&=&{6x\over \pi^2\prod_{p|r}(1+{1\over p})}+O(\sqrt{x}\varphi(r)),\nonumber
\end{eqnarray}
where we used that
$${\varphi(r)\over r}\sum_{(d,r)=1}^{\infty}{\mu(d)\over d^2}={\varphi(r)\over r}
\prod_{p\nmid r}(1-{1\over p^2})={\varphi(r)\over \zeta(2)r\prod_{p|r}(1-{1\over  p^2})}
={1\over \zeta(2)\prod_{p|r}(1+{1\over p})}$$
and $\zeta(2)=\pi^2/6$. \qed\\

\noindent We also need to estimate $M_r(x)=\sum_{m\le x\atop (m,r)=1}\mu(m)$. Estimating this quantity turns out
to be much harder. The quantity $M_1(x)$ is the summatory function of the M\"obius function and this function
is called the {\it Mertens function}. The prime number theorem is known to be equivalent with the assertion that
$\lim_{x\rightarrow \infty}M_1(x)/x=0$. The celebrated Mertens' conjecture states that $|M_1(x)|
<\sqrt{x}$ for every $x$. If this conjecture would hold, then from the easily proved equality
$${1\over \zeta(s)}=s\int_1^{\infty}{M_1(x)\over x^{s+1}}ds,$$
it would follow that $\zeta(s)$ has no zeros with Re$(s)>1/2$, i.e. the Riemann Hypothesis would follow.
In fact, the Riemann Hypothesis is known to be equivalent with the assertion that $M_1(x)=O(x^{1/2+\epsilon}),$
for every $\epsilon>0$. In 1986 Odlyzko and te Riele disproved the Mertens conjecture. The behaviour
of $M_1(x)$ is thus closely related to the zero free region of the Riemann zeta function. The largest
known zero free region of $\zeta(s)$, due to Korobov, was used by Walfisz to show that
$$M_1(x)=O(x\exp(-C(\log x)^{3/5}(\log \log x)^{-1/5})).$$
In this work we will be satisfied with merely determining the average of $M_r(x)$:
\begin{Lem}
\label{mobiustwee}
We have $$\sum_{m\le x\atop (m,r)=1}\mu(m)=o(x),$$ where the implied constant may depend on $r$.
\end{Lem}
To prove the lemma, we will
apply the Wiener-Ikehara Tauberian theorem in the following form.
\begin{Thm}
\label{vertaubung}
Let $f(s)=\sum_{n=1}^{\infty}a_n/n^s$ be a Dirichlet series. Suppose there exists a Dirichlet series
$F(s)=\sum_{n=1}^{\infty}b_n/n^s$ with positive real coefficients such that\\
{\rm (a)} $|a_n|\le b_n$ for all $n$;\\
{\rm (b)} the series $F(s)$ converges for Re$(s)>1$;\\
{\rm (c)} the function $F(s)$ can be extended to a meromorphic function in the region Re$(s)\ge 1$ having
no poles except for a simple pole at $s=1$.\\
{\rm (d)} the function $f(s)$ can be extended to a meromorphic function in the region Re$(s)\ge 1$ having
no poles except possibly for a simple pole at $s=1$ with residue $r$.\\
Then
$$\sum_{n\le x}a_n=rx+o(x),~x\rightarrow \infty.$$
In particular, if $f(s)$ is holomorphic at $s=1$, then $r=0$ and
$\sum_{n\le x}a_n=o(x)$ as $x\rightarrow \infty$.
\end{Thm} 
{\it Proof of Lemma} \ref{mobiustwee}. We apply the Wiener-Ikehara theorem with $F(s)=\zeta(s)$ and
$$f(s)=\sum_{(n,r)=1}{\mu(n)\over n^s}={1\over \zeta(s)\prod_{p|r}(1-{1\over p^s})}.$$
Of course $F(s)$ satisfies the required properties and has a simple pole at $s=1$ with residue one.
Since the finite product in the formula for $f(s)$ is regular for Re$(s)>0$, the result follows
on using the well-known fact that $1/\zeta(s)$ can be extended to a meromorphic function in
the region Re$(s)\ge 1$ (and hence $r=0$). \qed

\subsection{Counting squarefree shifted primes}
\noindent In this section we are concerned with estimating the second sum in (\ref{twosums}).
\subsubsection{On a result of Mirsky}
Theorem \ref{mirsky} below is due to Mirsky \cite{M} 
(with the difference that in his version the $O$-constant
depends at most on $k,r$ and $H$). In his paper Mirsky states
two theorems, of which 
he only proves the first (the proof of the second being similar).
Mirsky's second theorem is stated below.
For completeness we give the
proof here.
Recall that the letters $p$ and $q$ are used to indicate primes.
\begin{Thm}
\label{mirsky}
Let $r$ be any non-zero integer, $k$ any integer greater than $1$, and
$H$ any positive number. Then
$$\#\{q\le x:q-r~{\rm is~}k{\rm -free}\}=
\prod_{p\nmid r}\left(1-{1\over p^{k-1}(p-1)}\right){\rm Li}(x)
+O\left({x\over \log^H x}\right),$$
where the $O$-constant depends at most on $k$ and $H$.
\end{Thm}
Our proof rests on the Siegel-Walfisz theorem:
\begin{Lem} 
\label{siegelwalfisz}
{\rm \cite[Satz 4.8.3]{Prachar}}. 
Let $\pi(x;m,l)$ denote the number of primes $q\le x$ with
$q\equiv l({\rm mod~}m)$.
Let $C>0$ be arbitrary. Then
$$\pi(x;m,l)={{\rm Li}(x)\over \varphi(m)}+O(xe^{-c_1\sqrt{\log x}}),$$
uniformly for $1\le m\le \log^C x$, $(l,m)=1$, where the constants
depend at most on $C$.
\end{Lem}

\noindent Note that if $(l,m)>1$ there is at most one prime $p\equiv l({\rm mod~}m)$. If $(l.m)=1$, then
the above result implies that $\delta(p\equiv l({\rm mod~}m))=1/\varphi(m)$, thus asymptotically the
primes are equidistributed over the primitive congruence classes mod $m$.\\

\noindent {\it Proof of Theorem} \ref{mirsky}. 
The dependence of O-constants on $r,k$ and $H$ is indicated
by using them as index, thus $O_H$ means that the implied constant
depends at most on $H$.\\
\indent Let  $y=\log^{H}x$.
By Proposition \ref{kmoebius} we have
\begin{eqnarray}
\label{nogeenmirsky}
\#\{q\le x:q-r{\rm ~is~}k{\rm -free}\}&=&\sum_{q\le x}
 \sum_{a^k|q-r}\mu(a)\cr
 &=&\sum_{a\le y}\mu(a)\sum_{q\le x\atop q\equiv r({\rm mod~}a^k)}1
 +\sum_{a> y}\mu(a)\sum_{q\le x\atop q\equiv r({\rm mod~}a^k)}1\cr
 &=& I_1+I_2,
\end{eqnarray}
say. We have, using Lemma \ref{siegelwalfisz} with $C=kH$, 
\begin{eqnarray}
\label{mirskynul}
I_1&=&\sum_{a\le y\atop (a,r)=1}
\mu(a)\sum_{q\le x\atop q\equiv r({\rm mod~}a^k)}1+O(y)\cr
&=&\sum_{a\le y\atop (a,r)=1}
\mu(a)\Big\{{{\rm Li}(x)\over 
\varphi(a^k)}+O_H({x\over \log^{2H} x})\Big\}+O(y)\cr
&=&\sum_{a\le y\atop (a,r)=1}{\mu(a)\over \varphi(a^k)}{\rm Li}(x)
+O_H({xy\over \log^{2H} x})+O(y).
\end{eqnarray}
Clearly
\begin{equation}
\label{mirskyeen}
\sum_{a\le y\atop (a,r)=1}{\mu(a)\over \varphi(a^k)}=
\sum_{a=1\atop (a,r)=1}^{\infty}{\mu(a)\over \varphi(a^k)}
+O\left(\sum_{a\ge y}{1\over \varphi(a^k)}\right).
\end{equation}
By Euler's product identity we have
\begin{equation}
\label{mirskytwee}
\sum_{a=1\atop (a,r)=1}^{\infty}{\mu(a)\over \varphi(a^k)}=
\prod_{p\nmid r}\left(1-{1\over p^{k-1}(p-1)}\right).
\end{equation}
Using 
that $k\ge 2$ and the classical estimate $\varphi(m)\gg m/\log \log m$, we obtain
\begin{equation}
\label{mirskydrie}
\sum_{a\ge y}{1\over \varphi(a^k)}=O\left(\sum_{a\ge y}{\log (k\log a )\over
a^k}\right)=O_k\left({\log \log y\over y^{k-1}}\right)
=O_k\left({\log \log y\over y}\right).
\end{equation}
On combining
(\ref{mirskynul}) with (\ref{mirskyeen}), (\ref{mirskytwee}) and
(\ref{mirskydrie}), we infer that
$$I_1=\prod_{p\nmid r}\left(1-{1\over p^{k-1}(p-1)}\right){\rm Li}(x)
+O_k\left({x\log \log y\over (\log x) y}\right)+O_H({xy\over \log^{2H}x})
+O(y).$$ 
Note that
$$|I_2|\le \sum_{a>y}\sum_{m\le x\atop m\equiv r({\rm mod~}a^k)}1
\le \sum_{a>y}[{x\over a^k}]=O({x\over y^{k-1}})=O({x\over y}),$$
where in the second sum the summation is over all integers $m\le x$
satisfying $m\equiv r({\rm mod~}a^k)$. On adding the estimates
for $I_1$ and $|I_2|$, the result follows from 
(\ref{nogeenmirsky}). \qed

\subsubsection{Connection with Artin's primitive root conjecture}
\label{return}
\noindent By Theorem \ref{mirsky} the density of primes $p$ such that $p-1$
is squarefree equals $A$, the Artin constant. The Artin constant also
arose in Section \ref{primitiverOOts} in the context of primitive roots. This raises
the question whether there is some connection between the two problems. We will
now show that this is indeed the
case; both problems are in fact special cases of a
generalization of Artin's primitive root conjecture. We assume some familiarity with algebraic number theory.\\
\indent We recall that a prime $p$ splits completely in $\mathbb Q(\zeta_n)$  iff $p\equiv 1({\rm mod~}n)$.
Note that $p-1$ is squarefree iff $p\not\equiv 1({\rm mod~}q^2)$ with $q$ any prime. It follows that
$p-1$ is squarefree iff $p$ does not split completely in any of the fields $\mathbb Q(\zeta_{q^2})$ with
$q$ a prime. It is a consequence of Chebotarev's density theorem that the density of primes $q$ that
split completely in a normal extension $K:\mathbb Q$ equals $1/[K:\mathbb Q]$, where $[K:\mathbb Q]$ denotes
the degree of the extension. On noting that the compositum 
of $\mathbb Q(\zeta_{q_1^2}),...,\mathbb Q(\zeta_{q_s^2})$, where $q_1,...,q_s$ are distinct primes, 
equals $\mathbb Q(\zeta_{(q_1q_2\cdots q_s)^2})$, we expect by inclusion and exclusion that the density
of primes $p$ for which $p-1$ is squarefree equals
$$\sum_{n=1}^{\infty}{\mu(n)\over [\mathbb Q(\zeta_{n^2}):\mathbb Q]}
~(=\sum_{n=1}^{\infty}{\mu(n)\over \varphi(n^2)}=A).$$
Indeed, the primes $p$ for which a given integer $g$ with $g\ne -1$ or a square is a primitive
root mod $p$ can be described
in a similar way. Here we want that for each prime $q$ with $q|p-1$ we have that 
$g^{p-1\over q}\not\equiv 1({\rm mod~}p)$, which can be reformulated as $p$ does not split
completely in $\mathbb Q(\zeta_q,g^{1/q})$ for any
prime $q$. By inclusion and exclusion we expect then that that
the density of primes $p$ for which $g$ is a primitive root equals
\begin{equation}
\label{hooley2}
\sum_{n=1}^{\infty}{\mu(n)\over [\mathbb Q(\zeta_n,g^{1/n}):\mathbb Q]}.
\end{equation}
In fact Hooley proved in 1967 that the number given in
(\ref{hooley2}) is indeed the correct density, under the assumption of GRH.
For squarefree $n$ the degree $[\mathbb Q(\zeta_n,g^{1/n}):\mathbb Q]$ is generically equal to $n\varphi(n)$, 
which is the exact degree of $[\mathbb Q(\zeta_{n^2}):\mathbb Q]$ and hence the arisal of $A$ in
both problems does not surprise us.\\
\indent Both problems considered here are in fact special cases of the following generalization of
Artin's primitive root problem. This generalization was first studied by Goldstein and later
by M. Ram Murty \cite{Murty2}. Let $K$ be an algebraic number field. Let $\cal F$ be a family
of number fields, normal and of finite degree over $K$. Determine the number of prime ideals $\cal P$
of $K$ such that $N_{K/\mathbb Q}({\cal P})\le x$ and which do not split completely in any element
$\ne K$ of $\cal F$. On taking $K=\mathbb Q$ and $F$ to be the set of fields of the form
$\mathbb Q(\zeta_q,g^{1/q})$, where $q$ runs over the primes, we obtain Artin's primitive root problem.
On taking   $K=\mathbb Q$ and $F$ to be the set of fields of the form
$\mathbb Q(\zeta_{q^2})$, where $q$ runs over the primes, we obtain the special case
of Mirsky's result considered above.\\
\indent From an heuristical point of view the arisal of Artin's constant in the squarefree problem is not
surprising. The density of primes $q$ such that $q\not\equiv 1({\rm mod~}p^2)$ is $1-1/(p(p-1))$. Imposing
the condition $q\not\equiv 1({\rm mod~}p^2)$ for each prime $p$ and assuming all the conditions to represent independent events, one
arrives at $A$ as the expected density of the primes $q$ with $q-1$ squarefree. The independence of the
various local conditions is a consequence of the fact that for any two distinct primes $r$ and $s$ we have
$\mathbb Q(\zeta_{r^2}) \cap \mathbb Q(\zeta_{s^2})=\mathbb Q$.\\
\indent For the Artin problem we impose the local condition that $q$ does not satisfy both
$q\equiv 1({\rm mod~}p)$ and $g^{(q-1)/p}\equiv 1({\rm mod~}p)$. The density of primes $q$ not satisfying
the latter two conditions is $1-1/[\mathbb Q(\zeta_p,g^{1/p}):\mathbb Q]$. If the conditions for the
various primes $p$ would be independent, we would expect (as Artin original did) the density to equal
\begin{equation}
\label{verwachting}
\prod_p\left(1-{1\over [\mathbb Q(\zeta_p,g^{1/p}):\mathbb Q]}\right),
\end{equation}
which if $g\ne -1,0,1$ equals a rational number times the Artin constant.
However, it is not always true (with $r$ and $s$  as before) that $\mathbb Q(\zeta_r,g^{1/r})\cap \mathbb Q(\zeta_s,g^{1/s})=\mathbb Q$ 
(for example when $r=5$, $s=2$ and $g=5$ in which case the intersection equals $\mathbb Q(\sqrt{5})$) and
hence we expect the density to come out different from (\ref{verwachting}) sometimes (as can be shown
true assuming GRH).

\subsubsection{A variation of Mirsky's result}
We will need a variation of Mirsky's result for $k=2$ and $r=1$. Let $S=(q_1,\dots,q_t)$ be a (possibly empty)
finite sequence of primes satisfying $q_1<q_2< \dots < q_t$. We write $\nu_S(n)=(\nu_{q_1}(n),\dots,\nu_{q_t}(n))$. 
We say that a number $n$ is $S$-squarefree if $q^2\nmid n$ for every prime $q$ with $q\not\in S$. We put
$$\mu_S(n)=\cases{(-1)^{\sum_{q|n,~q\not\in S}1} &if $n$ is $S$-squarefree;\cr
0 & otherwise.}$$
If $S$ is the empty set, then $\mu_S(n)=\mu(n)$.
\begin{Prop}
\label{mirsky2}
Let $t$ be any non-zero integer, Let $S=\{q_1,\dots,q_t\}$ and $e_1,\dots,e_t$ be natural numbers. 
Write $Q=\prod_{j=1}^{t}q_j^{e_j}$.
We have
$$\sum_{p\le x\atop \nu_S(p-1)=(e_1,\dots,e_t)}\mu_S(p-1)^2=$$
$${1\over Q}\prod_{1\le j\le t\atop e_j=0}\left(1-{1\over 
q_j-1}\right)
\prod_{q\not\in S}\left(1-{1\over q(q-1)}\right){\rm Li}(x)
+O\left({x\over \log^H x}\right),$$
where the $O$-constant depends at most on $S,e_1,\dots,e_t$ and $H$.
\end{Prop}
{\it Proof}. A variation of Theorem \ref{mirsky} and left to the reader. \qed\\

\noindent {\tt Remark}. Again the constant appearing in this result is not surprising. For each prime $q$ we impose a condition
at our prime $p$ occurring in the sum: if $q\in S$ we prescribe $\nu_q(p-1)$ and if $q\not\in S$ we require
that $p\not\equiv 1({\rm mod~}q^2)$. 
If for fixed $q$ we compute the densities of primes satisfying this condition
and if we multiply all these `local densities' together, we arrive at the density
$$\delta(\mu_S(p-1)\ne 0)\prod_{j=1}^t \delta(\nu_{q_j}(p-1)=e_j)=$$
$${1\over \varphi(Q)}\prod_{1\le j\le t\atop e_j=0}\left(1-{1\over q_j-1}\right)
\prod_{1\le j\le t\atop e_j\ge 1}\left(1-{1\over q_j}\right)
\prod_{q\not\in S}\left(1-{1\over q(q-1)}\right),$$
which equals the density given in Proposition \ref{mirsky2}

\subsubsection{A conjecture}
It is generally believed that primes $p$ are such that $p-1$ does not have a preference with regards to having
an even or odd number of prime factors. In other words, the following conjecture (thought
to be deep by the experts) is generally believed:
we have
$\sum_{p\le x}\mu(p-1)=o(\pi(x)),$
as $x$ tends to infinity. We propose a slightly more general conjecture.
\begin{Con}
Let $t$ be any non-zero integer, Let $S=\{q_1,\dots,q_t\}$ and $e_1,\dots,e_t$ be natural numbers. 
We have
$$\sum_{p\le x\atop \nu_S(p-1)=(e_1,\dots,e_t)}\mu_S(p-1)=o(\pi(x)),$$
where the $o$-constant depends at most on $S,e_1,\dots,e_t$ and $H$.
\end{Con}
All the numerical data we came across in this respect seemed to be not inconsistent with this conjecture.

\section{Value distribution of Ramanujan sums}
\label{aurelius}
We discuss the value distribution of $c_n(m)$ as $n$ varies over the integers and $m$ is fixed. We
do this following Aurel Wintner's paper \cite{W}. Wintner argues that the prime number theorem is
equivalent with the M\"obius function having an asymptotic distribution function and notes that
$c_n(1)=\mu(n)$. A natural question arising then is whether $c_n(m)$ for $m$ fixed has
an asymptotic distribution function and if so, what it looks like. We next recall some facts 
regarding moments and distribution functions.

\subsection{Distribution functions}
A sequence $\{\alpha_j\}_{j=0}^{\infty}$ is said to have an {\it asymptotic distribution function}, $\rho$, 
if there exists a monotone function $\rho(\alpha)$ such that 
$$\lim_{x\rightarrow \infty}{\sum_{\alpha_j\le x}1\over x}=\rho(\alpha),$$
holds at every continuity point $\alpha$ of $\rho$ and, moreover, we have $\rho(-\infty)=0$ and
$\rho(\infty)=1$. It is known that if $\lim_{j\rightarrow \infty} \sup |\alpha_j|<\infty$, it has
an asymptotic distribution $\rho$ iff $M_n(\alpha_n^j)$ exists for every $j$; in which case
$$M_n(\alpha_n^j)=[\rho]_j=\int_{-\infty}^{\infty}\alpha^jd\rho(\alpha).$$
(We define $\alpha_n^0=1$ and hence $[\rho]_0=1$.)
Given a  distribution function $\rho(\alpha)$ we can consider, $L(u;\rho)$, the Fourier-Stieltjes
transform of $\rho$:
$$L(u;\rho)=\int_{-\infty}^{\infty}e^{i \alpha u}d\rho(\alpha).$$
If the moments exist for every $j$, then
$$L(u;\rho)=\sum_{j=0}^{\infty}{i^j[\rho]_j\over j!}u^j$$
is valid for every $u$. In the case where there are only finitely many distinct elements in ths sequence 
(which will happen in our application), 
$\rho(\alpha)$ will be a step function making only finitely many steps. Then $L(u;\rho)$ will be a finite
expression of the form $\sum_{j=1}^r c(\beta_j)e^{i\beta_r}$ and the values assumed are then
$\beta_1,\cdots,\beta_r$, each with density $c(\beta_j)$.\\
\indent As an illustration let us compute the asymptotic distribution for $\{\mu(n)\}_{n=1}^{\infty}$.
By now the reader should have no difficulties in proving directly that this is merely the step function which
has the saltus (`jump') $3/\pi^2,~1-6/\pi^2,3/\pi^2$ at $\alpha=-1,0,1$, respectively. Let us derive this
by the method of moments. The odd moments of $\mu$ are zero and the $2j$th moment with $j\ge 1$ equals
$6/\pi^2$. The zeroth moment equals one. We thus derive that 
$$L(u,\rho)=1-{6\over \pi^2}+{6\over \pi^2}\cos u={3\over \pi^2}e^{-iu}+1-{6\over \pi^2}+{3\over \pi^2}e^{iu},$$
from which we arrive at the same conclusion as before.

\subsection{Computing the moments and asymptotic distribution function}
Using the representation (\ref{wrelatie}) for the Dirichlet series $f_n^{(2j)}(s)$ and
$f_n^{(2j-1)}(s)$ and the Wiener Ikehara theorem, Theorem 
\ref{vertaubung}, the mean of the $j$th moment of $c_n(m)$
can be determined (note that by property 5 of Proposition \ref{basicramanujan} 
we can take $F(s)=m^j\zeta(s)$ in that result when we apply it to $f_n^{(j)}(s)$), We
find that
\begin{equation}
\label{nogeenwinter} 
M_n(c_n(m)^{2j})={F_m^{(2j)}(1)\over \zeta(2)}={6\over \pi^2}\prod_{p|m}F_{p^{\nu_p(m)}}^{(2j)}(1).
\end{equation}
Note that $M_n(c_n(m)^{2j})\zeta(2)$ is a multiplicative function in $m$.
Similarly, we have $M_n(c_n(m)^{2j-1})=0$. (We use the subscript $n$ in $M$ to indicate that the mean is 
with respect to $n$.) Since $|c_n(m)|^{2j}\le m^{2j}$ and hence the same inequality holds for the mean,
it follows that the radius of convergence of the power series
\begin{equation}
G_m(u)=\sum_{j=0}^{\infty}{M_n(c_n(m)^{2j})\over (2j)!}(-u^2)^j
\end{equation}
is infinite for every $m~(=1,2,\cdots)$. In case $m=p^r$, $G_m(u)$ is easily evaluated using
(\ref{nogeenwinter}) and (\ref{eerstewinter}). One obtains
$$G_{p^r}(u)=\left(1-{6\over \pi^2}{1-p^{-k-2}\over 1-p^{-2}}\right)+{6\over \pi^2}\left(1+{1\over p}\right)^{-1}
\left(\sum_{j=0}^r {\cos (\varphi(p^j)u)\over p^j}+{\cos(p^ru)\over p^{r+1}}\right).$$
Note that if $k$ is fixed, then
$$G_{p^r}(u)\rightarrow 1-{6\over \pi^2}+{6\over \pi^2}\cos u,~(p\rightarrow \infty)=G_1(u),$$
and
$$G_{p^r}(u)\rightarrow \left(1-{6\over \pi^2(1-p^{-2})}\right)+
{6\over \pi^2}(1+{1\over p})^{-1}
\sum_{j=0}^r {\cos (\varphi(p^j)u)\over p^j},~(k\rightarrow \infty).$$
For any given $m$ it is not difficult to compute $G_m(u)$ using formula (\ref{nogeenwinter}).
The existence of all the moments $M_n(c_n(m)^j)$, $j=0,1,2,\cdots$ implies the existence of 
an asymptotic distribution function $\rho_m$. Note that $L(u;\rho_m)=G_m(u)$ and that
$G_m(u)=\sum_{j=1}^r c(\beta_j)e^{i\beta_r}$. The values assumed are 
$\beta_1,\cdots,\beta_r$, each with density $c(\beta_j)$. Thus, taking $m=p^r$ for example we
immediately read off from the formula for $G_{p^r}(u)$ that 
$$
\delta(c_n(p^r)=v)=\cases{1-{6\over \pi^2}{1-p^{-r-2}\over 1-p^{-2}} &if $v=0$;\cr
{3\over \pi^2 p^h(1+{1\over p})} &if $v=\pm \varphi(p^h)$, where $h=1,\dots,r$;\cr
{3\over \pi^2 p^{r+1}(1+{1\over p})} &if $v=\pm p^r$.}$$


\section{Value distribution of cyclotomic coefficients}
\label{vdofcyclotomic}
\subsection{Evaluating cyclotomic coefficients}
Write $$\Phi_n(X)=\sum_{k=0}^{\varphi(n)}a_n(k)X^k.$$ 
We consider the value distribution of the cyclotomic coefficients $a_n(k)$ for fixed $k$, as
$n$ runs over all the positive integers. 
By setting $\mu(n/d)=0$ whenever $n/d$ is not an integer
we have, for $n>1$,
\begin{equation}
\label{thanga}
\Phi_n(X)=\prod_{d|n}(1-X^d)^{\mu(n/d)}=\prod_{d=1}^{\infty}(1-X^d)^{\mu(n/d)}.
\end{equation}
Note that, for $|X|<1$, we have
\begin{equation}
\label{lehmertje}
\prod_{d=1}^{\infty}(1-X^d)^{\mu(n/d)}=\prod_{d=1}^{\infty}\left(1-\mu({n\over d})X^d+
{1\over 2}\mu({n\over d})(\mu({n\over d})-1)\sum_{j=1}^{\infty}X^{jd}\right),
\end{equation}
where we used the observation that, for $|X|<1$,
\begin{equation}
\label{thanga2}
(1-X^d)^{\mu(n/d)}=\left(1-\mu({n\over d})X^d+
{1\over 2}\mu({n\over d})(\mu({n\over d})-1)\sum_{j=1}^{\infty}X^{jd}\right).
\end{equation}
From (\ref{lehmertje}) it is not difficult to derive a formula for $a_n(k)$ for a fixed $k$;
this is just the coefficient of $X^k$ in the right hand side of (\ref{thanga}) (this approach
seems to be due to D.H. Lehmer \cite{DL}). We thus obtain, 
for $n>1$, 
$$\cases{
a_n(1)=-\mu(n);\cr
a_n(2)=\mu(n)(\mu(n)-1)/2-\mu(n/2);\cr
a_n(3)=\mu(n)^2/2-\mu(n)/2+\mu(n/2)\mu(n)-\mu(n/3).}$$
More generally, we have
\begin{equation}
\label{uitdrukking}
a_n(k)=\sum c(k_1,...,k_s;e_1,\dots,e_s)\mu({n\over k_1})^{e_1}\cdots \mu({n\over k_s})^{e_s},
\end{equation}
where the sum is over all partitions $k_1+\cdots+k_s$ of all the
integers $\le k$ with $k_1\ge k_2\ge \cdots k_s$ and over all $e_1,\dots,e_s$ with
$1\le e_j\le 2$ for $1\le j\le s$.
Using this result we will deduce that the $n$ dependence of $a_n(k)$ is not that strong.
\begin{Prop}
\label{naareindig}
Put $N_k={\rm lcm}(1,2,\cdots,k)\prod_{p\le k}p$. We can uniquely decompose $n$
as $n=n_kc_k$ with $(c_k,N_k)=1$
and $n_k$ and $c_k$ natural numbers.
There exist functions $A_1$ and $B_1$ with as domain the divisors of $N_k$ such that
$$a_n(k)=\cases{A_1(n_k)\mu(c_k)^2+B_1(n_k)\mu(c_k) &if $n_k|N_k$;\cr
0 &otherwise.}$$
\end{Prop}
{\it Proof}. The assertion regarding the uniqueness of the decomposition $n=n_kc_k$ is trivial. 
For a given $n$ only those partitions $k_1,k_2,\cdots,k_s$ contribute to 
(\ref{uitdrukking}) for which
$n/k_i$ is an integer for $1\le i\le s$. Note that $k_i|n_k$. Thus, we can write
$$\mu({n\over k_1})^{e_1}\cdots \mu({n\over k_s})^{e_s}=
\mu({n_k\over k_1})^{e_1}\cdots \mu({n_k\over k_s})^{e_s}\mu(c_k)^{e_1+\dots+e_s}.$$
If $n_k\nmid N_k$, then none of the integers $n_k/k_1,...,n_k/k_s$ is squarefree and 
$a_n(k)=0$, so assume that $n_k|N_k$. 
On using that $\mu(r)^w$ with $w\ge 1$ either equals $\mu(r)$ or $\mu(r)^2$, the result follows
from (\ref{uitdrukking}). \qed\\


\noindent The above proposition shows that $\{a_n(k)|n\in \mathbb N\}$ is a finite set, thus if
we fix $k$, there are only finitely many possibilities for the values of the
coefficient of $X^k$ in a cyclotomic polynomial. We will now show that
$-1,0$ and $1$ are always amongst these values. In the formula for $a_n(k)$ there is always
the term $-\mu(n/k)$. Let us take $n=ck\prod_{p\le k}p$, where $c$ only has
prime divisor $>k$. Then all the terms of the form $\mu(n/r)$
with $1\le r<k$ are zero (since either $r\nmid n$ or $n/r$ is not squarefree) and we obtain
that $a_n(k)=-\mu(c)(-1)^{\pi(k)}$, where $\pi(x)$ as usual denotes the number of primes $p\le x$
not exceeding $x$. In particular, it follows that $a_n(k)$ always assumes the values $-1,0$ and $1$. 
Combining this insight with the previous proposition we arrive at the following result.
\begin{Prop}
\label{valueset}
Let $k\ge 1$ be fixed. Let $N_k={\rm lcm}(1,2,\cdots,k)\prod_{p\le k}p$ and 
let $q>k$ be any prime exceeding $k$. We have
$$\{-1,0,1\}\subseteq \{0,a_{d}(k),a_{dq}(k)~|~d|N_k\}=\{a_n(k)~|~n\ge 1\}.$$
\end{Prop}
Using formula (\ref{naarkwadraatvrijekern}) it is seen that in the latter proposition $N_k$ can
be replaced by $k\prod_{p\le k}p$.
We have $|a_n(k)|\le \max_{n\ge 1}|a_n(k)|=B(k)$. See Table 2 for the values
of $B(k)$ for $1\le k\le 30$.\\
\indent We next show that the inclusion in Proposition \ref{valueset} is strict for $k\ge 13$. Our
proof rests on 
the following rather elementary result on prime numbers. 
\begin{Prop}
\label{bertie}
For $k\ge 13$ there are consecutive odd primes $p_1<p_2<p_3$ such that $p_3\le k<p_1+p_2$.
\end{Prop}
{\it Proof}. Breusch proved in 1934 that for $x\ge 48$ there is at least one prime in $[x,9x/8]$ 
(this strengthens {\it Bertrand's Postulate} asserting that there is always a prime between $x$ and
$2x$, provided $x\ge 2$).
Let $\alpha=1.32$. A little computation shows that the above result implies that for $x\ge 9$ there
is at least one prime in $[x,\alpha x]$. One checks that the assertion is true for $k\in [13,21)$.
Assume that $k\ge 21~(\ge 9\alpha^3)$. Let $p_3$ be the largest prime not exceeding $k$ and let
$p_1$ and $p_2$ be primes such that $p_1,p_2$ and $p_3$ are consecutive primes. Then
$p_3\ge k/\alpha$, $p_2\ge k/\alpha^2$ and $p_1\ge k/\alpha^3$. On noting that
$p_1+p_2\ge k(1/\alpha+1/\alpha^2)>k$, the proof is then completed. \qed

\begin{Prop}
For $k\ge 13$ we have $\{-2,-1,0,1\}\in \{a_n(k):n\in \mathbb N\}$ (and thus $B(k)\ge 2$).
\end{Prop}
{\it Proof}. Let $p_1,p_2$ and $p_3$ be odd primes satisfying the condition of Proposition \ref{bertie}. 
Using  (\ref{thanga}) we infer that
$$\Psi_{p_1p_2p_3}(X)\equiv {(1-X^{p_3})\over (1-X)}(1-X^{p_1})(1-X^{p_2})\equiv (1+X+\cdots+X^{p_3-1})
(1-X^{p_1}-X^{p_2}),$$
where we computed modulo $X^{k+1}$.
It follows that $a_{p_1p_2p_3}(k)=-2$. The proof is completed on invoking
Proposition \ref{valueset}. \qed\\

\subsubsection{Numerical evaluation of $a_n(k)$ for small $k$}
For our purposes it is relevant to be able to 
numerically evaluate $a_n(k)$ for small $k$ and large $n$. A computer package like Maple
evaluates $a_n(k)$ be evaluating the whole polynomial $\Phi_n(x)$. For large $n$ this
is far too costly. Instead it is more efficient to use (\ref{thanga}) and expand for
every $d$ for which $\mu(n/d)\ne 0$, $(1-X^d)^{\mu(n/d)}$ as a Taylor series up to $O(x^{k+1})$
and multiply all these series together. The most efficient method to date is due to
Grytczuk and Tropak \cite{GT}. First they apply formula (\ref{naarkwadraatvrijekern}). Thus
it is enough to compute $a_n(k)$ with $n$ squarefree. If $\phi(n)<k$, then $a_n(k)=0$, so
we may assume that $\phi(n)\ge k$. Let $d=(n,\prod_{p\le k}p)$. Put
$T_r=\mu(n)\mu((r,d))\varphi((r,d))$. Compute $b_0,\dots,b_k$ recursively by $b_0=1$ and
$$b_j=-{1\over j}\sum_{m=0}^{j-1}b_mT_{j-m}{\rm ~for~}j>0.$$
Then $b_k=a_n(k)$. The proof uses the formula
\begin{equation}
\label{tropakje}
a_n(k)=-{1\over k}\sum_{m=0}^{k-1}a_n(m)c_n(k-m)~{\rm for~}k>1,
\end{equation}
which follows by Viete's and Newton's formulae from (\ref{definitie}) and it uses Property 2
of Proposition \ref{basicramanujan}. An alternative proof of (\ref{tropakje}) is obtained
on using the following observation together with  Property 1
of Proposition \ref{basicramanujan}.
\begin{Prop}
Suppose that, as formal power series, 
$$\prod_{d=1}^{\infty}(1-X^d)^{-a_d}=\sum_{d=0}^{\infty}r(d)X^d,$$
then $dr(d)=\sum_{j=1}^d r(d-j)\sum_{k|j}ka_k$.
\end{Prop}
{\it Proof}. Taking the logarithmic derivative of $\prod_{d=1}^{\infty}(1-X^d)^{-a_d}$ we obtain
$${\sum_{d=1}^{\infty}dr(d)X^d\over \sum_{d=0}^{\infty}r(d)X^d}=X{d\over dx}\log 
\prod_{d=1}^{\infty}(1-X^d)^{-a_d}=\sum_{j=1}^{\infty}(\sum_{k|j}ka_k)X^j,$$
whence the result follows. \qed\\

\noindent For every integer $v$ it is a consequence of Proposition 
\ref{z} that there exists a minimal integer
$k$, $k_{\rm min}$, such that there exists a natural number $n$ with $a_n(k_{\rm min})=v$.  
Grytczuk and Tropak \cite[Table 2.1]{GT} used their method to determine $k_{\rm min}$ for the integers in the
interval $[-9,\dots,10]$. Gallot has extended this range from $[-60,\dots,70]$. From Table 11
one can determine $k_{\rm min}$ for the range $[-15,\dots,15]$.
\subsection{M\"oller's paper reconsidered}
In this section we reconsider M\"oller's paper \cite{HM}. 
In M\"oller's approach $a_n(k)$ is 
connected with partitions of $k$. A partition of a positive integer $m$
is an expression of the form $m_1+m_2+\cdots+m_r=m$ with all the $m_j\ge 0$. Ordering is
disregarded. Thus, $1+1+2+3$ and $3+2+1+1$ are considered to be the 
same partitions of 7. A partition can be
identified with a sequence $\{n_j\}_{j=0}^{\infty}$ of non-negative integers satisfying
$\sum_{j}jn_j=m$. 
W.l.o.g. we can denote a partition, $\lambda$, of $k$ as
$\lambda=(k_1^{n_{k_1}}\cdots k_s^{n_{k_s}})$, where
$n_{k_1}\ge n_{k_2}\ge \dots \ge n_{k_s}\ge 1$ (thus the number $k_j$ occurs $n_{k_j}$ times
in the partition). The set of all partitions of $m$ will be denoted by ${\cal P}(m)$. 
The number of different partitions of $m$ is denoted by $p(m)$. Hardy and
Littlewood in 1918, and Uspensky independently in 1920, proved that
$$p(m)\sim e^{\pi \sqrt{2m/3}}/(4m\sqrt{3})~{\rm ~as~}m\rightarrow \infty.$$
\begin{Prop} {\rm \cite[Satz 2]{HM}}.
\label{simpelzeg}
We have
$$a_n(k)=\sum_{\sum_{j}{jn_j}=k,~n_j\ge 0}\prod (-)^{n_j}\left({\mu({n\over j})\atop n_j}\right).$$
\end{Prop}
{\it Proof}. The taylor series of $(1-X)^a$ equals, for $|X|<1$, 
$\sum_{j=0}^{\infty}(-1)^j({a\atop j})X^j$, where
$({a\atop j})=a(a-1)\cdots (a-(j-1))/j!$. Using this we infer that
\begin{equation}
\label{nogeentje}
(1-X^d)^{\mu({n\over d})}=\sum_{j=0}^{\infty}(-1)^j\left({\mu({n\over d})\atop j}\right)X^{Dj},~|X|<1,
\end{equation}
The proof now follows from (\ref{nogeentje}) and (\ref{thanga}). \qed\\

\noindent Our proof above shows that the formula given in 
Proposition \ref{simpelzeg}  is a triviality, whereas M\"oller's 
ingenious and rather involved proof of it obscures this.\\
\indent M\"oller uses Proposition \ref{simpelzeg} to show that
$$M_n(a_n(k))=\lim_{x\rightarrow \infty}{\sum_{n\le x}a_n(k)\over x}$$
exists and gives a formula for it. To do so
he first 
determines the average of $\prod (-1)^{n_j}({\mu(n/j)\atop n_j})$. He does this by expanding it out
as a sum of terms of the form $\prod \mu^{\alpha_j}(n/j)$ with $1\le \alpha_j\le 2$, which is possible by
(\ref{veertien}). The average of each term $\prod \mu^{\alpha_j}(n/j)$ is easily determined. He then
invokes Proposition \ref{simpelzeg} and uses combinatorics to simplify his expressions. His final result is
still quite complicated and we did not feel inclined to check its equivalence with our more
simple formula in Theorem \ref{eenvoudig} below. We will now follow M\"oller's trail, in spirit,
but not in detail and see where it leads us.
\begin{Prop}
\label{halfweg}
Let $s\ge 1$, $k_1,\cdots,k_s$ be distinct integers and $n_{k_1}\ge n_{k_2}\ge \dots \ge n_{k_s}\ge 1$.
If $n_{k_1}\ge 2$ we let $t$ be the largest integer $\le s$ for which $n_{k_t}\ge 2$, otherwise
we let $t=0$.
Let $L=[k_1,\dots,k_s]$ and $G=(k_1,\dots,k_s)$. We have
$$\lim_{x\rightarrow \infty}{1\over x}\sum_{n\le x}(-1)^{n_{k_1}+\cdots+n_{k_s}}\left({\mu({n\over k_1})\atop n_{k_1}}\right)\cdots 
\left({\mu({n\over k_s})\atop n_{k_s}}\right)={6\over \pi^2}{\epsilon
\mu({L\over k_{t+1}})\cdots \mu({L\over k_s})
\over G\prod_{p|{L\over G}}(p+1)},$$
where
$$\epsilon=\cases{1 &if $n_{k_1}=1$, $s$ is even and $\mu(L/G)\ne 0$;\cr
\mu({L\over k_1})^{s-t}/2 &if $n_{k_1}\ge 2$ and $\mu(L/k_1)=\cdots=\mu(L/k_t)$ and $\mu(L/G)\ne 0$;\cr
0 & otherwise.}$$
\end{Prop}
\begin{cor}
\label{coreleganter}
We have $$\lim_{x\rightarrow \infty}{1\over x}\sum_{n\le x}(-1)^{n_{k_1}+\cdots+n_{k_s}}\left({\mu({n\over k_1})\atop n_{k_1}}\right)\cdots 
\left({\mu({n\over k_1})\atop n_{k_1}}\right)=$$
$${3\over \pi^2 G\prod_{p|{L\over G}}(p+1)}\left(\prod_{j=1}^s (-1)^{n_{k_j}} ({\mu({L\over k_j})\atop n_{k_j}})
+  \prod_{j=1}^s (-1)^{n_{k_j}} ({-\mu({L\over k_j})\atop n_{k_j}})\right).$$
\end{cor}
\noindent {\it Proof of Proposition } \ref{halfweg}. 
Put 
$$S(x)=\sum_{n\le x}(-1)^{n_{k_1}+\cdots+n_{k_s}}\left({\mu({n\over k_1})\atop n_{k_1}}\right)\cdots 
\left({\mu({n\over k_s})\atop n_{k_s}}\right).$$
Comparison of (\ref{thanga2}) and (\ref{nogeentje}) yields
\begin{equation}
\label{veertien}
(-1)^j\left({\mu({n\over d})\atop j}\right)=\cases{1 & if $j=0$;\cr
-\mu(n/d) &if $j=1$;\cr
\mu(n/d)(\mu(n/d)-1)/2 &if $j\ge 2$.}
\end{equation}
Note that for $j\ge 2$ the binomial coefficient is only non-zero if $\mu(n/d)=-1$. 
Using (\ref{veertien}) it
follows that a necessery condition for the argument of $S(x)$ to be non-zero is that
$L|n$. Now write $n=mL$. Note that $\mu(mL/k_1)\cdots \mu(mL/k_s)=\mu(m)^s \mu(L/k_1)\cdots \mu(L/k_s)$
if $(m,L/k_j)=1$ for $1\le j\le s$ and equals zero otherwise. It is not difficult to show that
$[{L\over k_1},\dots,{L\over k_s}]={L\over G}$ and using this, that $\mu(L/k_1),\cdots,\mu(L/k_s)$ are
squarefree iff $L/G$ is  squarefree. It follows
that if $\mu(L/G)=0$, then $S(x)=0$ and we are done, so next assume
that $\mu(L/G)\ne 0$. We infer that
$$S(x)=\sum_{m\le x/L,~(m,L/G)=1}(-1)^{n_{k_1}+\cdots+n_{k_s}}\left({\mu(mL/k_1)\atop n_{k_1}}\right)\cdots 
\left({\mu(mL/k_s)\atop n_{k_s}}\right).$$
Let us first consider the (easy) case where $n_{k_1}\ge 1$. Then we obtain
$S(x)=(-1)^s\mu({L\over k_1})\cdots \mu({L\over k_s})\sum_{m\le x/L,~(m,L/G)=1}\mu(m)^s$.
If $s$ is odd, then by Lemma \ref{mobiustwee} it follows that $\lim_{x\rightarrow \infty}S(x)/x=0$
and we are done, so next assume that $s$ is even. Then we apply Lemma \ref{mobiuseen} and
obtain that
$$\lim_{x\rightarrow \infty}{S(x)\over x}={6 \mu({L\over k_1})\cdots \mu({L\over k_s})\over
\pi^2 L\prod_{p|{L\over G}}(1+{1\over p})}.$$
The assumption $\mu(L/G)\ne 0$ implies that  $L\prod_{p|L/G}(1+1/p)=G\prod_{p|L/G}(p+1)$.\\
Next we consider the case where $n_{k_1}\ge 1$. The corresponding binomial coefficient is
only non-zero if $\mu(mL/k_1)=-1$. Similarly, we must have
$\mu(mL/k_j)=-1$ for $1\le j\le t$. It follows that if it is not true that
$\mu(L/k_1)=\dots=\mu(L/k_t)$, then $S(x)=0$ and hence $\lim_{x\rightarrow \infty}S(x)/x=0$ 
as asserted, so assume that $\mu(L/k_1)=\dots=\mu(L/k_t)$. We have
\begin{eqnarray}
S(x)&=&\sum_{m\le x/L,~(m,L/G)=1\atop \mu(mL/k_1)=-1}(-\mu(m))^{s-t}
\mu({L\over k_1})\cdots \mu({L\over k_1})\nonumber\\
&=&(\mu({L\over k_1}))^{s-t}\mu({L\over k_{t+1}})\cdots \mu({L\over k_s})
\sum_{m\le x/L,~(m,L/G)=1\atop \mu(m)=-\mu(L/k_1)}1.\nonumber
\end{eqnarray}
On invoking Lemma \ref{mobiuseen} and Lemma \ref{mobiustwee} the proof
is then completed. \qed\\

\noindent {\it Proof of Corollary} \ref{coreleganter}. Follows from a case by case analysis from
Proposition \ref{halfweg} on using (\ref{veertien}). \qed

\begin{Thm}
\label{eenvoudig}
We have
$$M_n(a_n(k))={3\over \pi^2}\sum_{\lambda=(k_1^{n_{k_1}}\dots k_s^{n_{k_s}})\in {\cal P}(k)\atop
n_{k_1}\ge \dots n_{k_s}\ge 1}{\epsilon(\lambda)\over G(\lambda)\prod_{p|{L(\lambda)\over G(\lambda)}}(p+1)},$$
where
$$\epsilon(\lambda)=\prod_{j=1}^s (-1)^{n_{k_j}} \left({\mu({L\over k_j})\atop n_{k_j}}\right)
+  \prod_{j=1}^s (-1)^{n_{k_j}} \left({-\mu({L\over k_j})\atop n_{k_j}}\right),$$
$L(\lambda)=[k_1,\dots,k_s]$ and $G(\lambda)=(k_1,\dots,k_s)$.
\end{Thm}
{\it Proof}. The result follows from Proposition \ref{simpelzeg} together with
Corollary \ref{coreleganter}. \qed\\

\noindent In case $k$ is a prime, $G(\lambda)=1$ for every partition and the above formula
further simplifies. The above formula suggests a connection with the group or representation
theory of the symmetric group $S_k$. The conjugacy classes in $S_k$ are in 1-1 correspondence
with the partitions of $k$. If $\lambda=(k_1^{n_{k_1}}\dots k_s^{n_{k_s}}$, then the order
of every element in the corresponding conjugacy class equals $L(\lambda)$. In particular,
$L(\lambda)\le g(n)$, where $g(n)$ denotes the maximum of all orders of elements in $S_k$.
It was shown by E. Landau in 1903 that $\log g(k)\sim \sqrt{k\log k}$ as $k$ tends to infinity
(for a nice account of this see \cite{Miller}), whereas by Stirling's theorem 
$\log k! \sim k\log k$.\\
\indent From Proposition \ref{simpelzeg} M\"oller infers that $|a_n(k)|\le p(k)-p(k-2)$.
To see this note that the partitions having 1 occurring at least twice do not contribute if
$\mu(n)\in \{0,1\}$. If $\mu(n)=-1$, then either $\mu(2n)=1$ or $\mu(n/2)=1$. This in
combination with (\ref{verdubbeling}) allows us then to argue as before and leads
us to the same bound. By a much more involved argument, using some
analytic number theory, M\"oller concludes that $B(k)>k^m$ 
for $k\ge k_0(m)$. This result was sharpened by several later authors culminating
in Bachman's estimate (\ref{gennady}).\\
\indent M\"oller used Proposition \ref{simpelzeg} to evaluate $B(k)$ for $1\le k\le 30$. 
He did this by hand. The
outcome is in Table 2. We redid this computation by computer and arrived the same result.  
Note that $M_n(a_n(k))=6e_k/\pi^2$, with $e_k$  a rational 
number.  For $1\le k\le 20$ we give the value of $e_k$ in Table 3 (our table agrees with
the one given in \cite{HM}, except for the incorrect values $e_{10}=319/1440$ 
and $e_{16}=733/2016$ appearing there). \\

\centerline{{\bf Table 3:} Scaled average, $e_k=\zeta(2)M_n(a_n(k))$, of $a_n(k)$ }
\begin{center}
\begin{tabular}{|c|c|c|c|c|c|c|c|c|c|c|c|}
\hline
$j$  & $1$ & $2$ & $3$ & $4$ & $5$ & $6$ & $7$ & $8$ & $9$ & $10$ \\
\hline
$e_j$  & $0$ & ${1\over 2}$ & ${1\over 6}$ & ${1\over 3}$ & ${1\over 8}$ & 
${7\over 24}$ & ${1\over 18}$ & ${7\over 24}$ & ${19\over 144}$ & ${31\over 160}$\\
\hline
$j$  & $11$ & $12$ & $13$ & $14$ & $15$ & $16$ & $17$ & $18$ & $19$ & $20$\\
\hline
$e_j$  & ${1\over 16}$ & ${55\over 192}$ & ${13\over 288}$ & ${61\over 288}$ & ${2287\over 20160}$ &
${733\over 4032}$ & ${667\over 8064}$ & ${79\over 336}$ & ${55\over 1344}$ & ${221\over 960}$  \\
\hline
\end{tabular}
\end{center}
\medskip
Regarding $e_k$
M\"oller proposed:
\begin{Con} {\rm \cite{HM}}. 
\label{vermoedentwee}
Let $k\ge 1$. Write $M_n(a_n(k))=6e_k/\pi^2$.\\
{\rm 1)} We have $0\le e_k\le 1/2$.\\
{\rm 2)} We have $(-1)^k(e_k-e_{k+1})>0$.
\end{Con}
M\"oller stated that with help of an IBM 7090 he wanted
to check his conjecture for further values of $k$. Had he carried this out, he would have
discovered that $(-1)^{34}(e_{34}-e_{35})=-18059/4626720<0$. Other
counterexamples occur at $k=35,45$ and $94$. Indeed, we would not be surprised if part 2 of
the Conjecture is violated for infinitely many $k$.\\
\indent On the other hand, part 1 of the Conjecture is true for $k\le 75$. The numbers
$e_k$ seem to be decreasing to zero and their size seems to be related to the number
of prime factors of $k$, the more prime factors the larger $e_k$ seems to be. In Table 11
(kindly computed by Yves Gallot) we give values for $e_k$ for $k$ up to $61$.

\subsection{Average and value distribution}

\noindent We give, using
Proposition \ref{naareindig}, a simpler formula for $M_n(a_n(k))$ involving $a_n(k)$
for a finite set of $n$. 

\begin{Thm}
\label{vier}
Let $k\ge 1$ be fixed. Put $M_k=k\prod_{p\le k}p,$ and let
$q>k$ be any prime. Then
$$M_n(a_n(k))={3\over \pi^2\prod_{p\le k}(1+{1\over p})}\sum_{d|M_k}{a_{d}(k)+a_{dq}(k)\over d}.$$
Furthermore, when $v\ne 0$,
$$\delta(a_n(k)=v)={3\over \pi^2\prod_{p\le k}(1+{1\over p})}(\sum_{d|M_k\atop a_{d}(k)=v}{1\over d}+
\sum_{d|M_k\atop a_{dq}(k)=v}{1\over d}).$$
\end{Thm} 
{\it Proof}. Let $N_k={\rm lcm}(1,2,\cdots,k)\prod_{p\le k}p$,
$r_1=\prod_{p\le k}p$ and $r_2=r_1N_k$. We have
\begin{eqnarray}
\sum_{n\le x}a_n(k)&=&\sum_{d|N_k}\sum_{n\le x\atop (n,r_2)=d}(A_1(d)\mu({n\over d})^2+B_1(d)\mu({n\over d}))\nonumber\\
&=&\sum_{d|N_k}\sum_{m\le x/d\atop (m,r_2/d)=1}(A_1(d)\mu(m)^2+B_1(d)\mu(m))\nonumber\\
&=&\sum_{d|N_k}\sum_{m\le x/d\atop (m,r_1)=1}(A_1(d)\mu(m)^2+B_1(d)\mu(m))\nonumber\\
&=&\sum_{d|N_k}A_1(d)\sum_{m\le x/d\atop (m,r_1)=1}\mu(m)^2+o_k(x),\nonumber
\end{eqnarray}
where we used Proposition \ref{naareindig} and Lemma \ref{mobiustwee}.
On invoking Lemma \ref{mobiuseen} we then obtain that 
$$\sum_{n\le x}a_n(k)={6\over \pi^2\prod_{p\le k}(1+{1\over p})}\sum_{d|N_k}{A_1(d)\over d}+o_k(x).$$
On noting that
$A_1(d)=(a_{d}(k)+a_{dq}(k))/2$ (and $B_1(d)=(a_{d}(k)-a_{dq}(k))/2$, but this is not needed),
the first formula follows, but with $N_k$ instead of $M_k$. 
On invoking formula (\ref{naarkwadraatvrijekern}) it is easily seen that 
$$\sum_{d|N_k}{a_d(k)+a_{dq}(k)\over k}=\sum_{d|M_k}{a_{d}(k)+a_{dq}(k)\over k}$$
and hence the first formula follows.
The proof of the second formula is similar and left to the reader. \qed\\

\noindent Using identity (\ref{verdubbeling}) we arrive at the following corollary to this theorem:
\begin{cor}
\label{corvier}
Let $k\ge 3$ be fixed and odd and $M_k$ and $q$ as in Theorem {\rm \ref{vier}}. Then
$$M_n(a_n(k))={1\over \pi^2\prod_{2<p\le k}(1+{1\over p})}\sum_{d|M_k/2}{a_{d}(k)+a_{dq}(k)\over d}.$$
Furthermore, when $v\ne 0$,
$$\delta(a_n(k)=v)={3\over \pi^2\prod_{p\le k}(1+{1\over p})}(\sum_{d|M_k/2\atop a_{d}(k)=v}{1\over d}+
\sum_{d|M_k/2\atop a_{d}(k)=-v}{1\over 2d}+
\sum_{d|M_k/2\atop a_{dq}(k)=v}{1\over d}+\sum_{d|M_k/2\atop a_{dq}(k)=-v}{1\over 2d}).$$
\end{cor}

\noindent As a corollary of Theorem \ref{vier} and Corollary \ref{corvier} we have:
\begin{Prop}
We have $e_k2k\prod_{p\le k}(p+1)\in \mathbb Z$. In case $k$ is odd we have
$e_kk\prod_{p\le k}(p+1)\in \mathbb Z$.
\end{Prop}
Yves Gallot observed that actually for $k\le 100$, we have $e_kk\prod_{p\le k}(p+1)\in \mathbb Z$.\\

\noindent In case $k$ is prime, the divisor sum in the previous corollary can be further reduced:
\begin{Prop} 
\label{priempiet}
Let $k\ge 3$ be a fixed prime. Put $N_k=\prod_{2<p<k}p,$ and let
$q>k$ be any prime. Then
$$M_n(a_n(k))={1\over \pi^2\prod_{2<p< k}(1+{1\over p})}\sum_{d|N_k}{a_{d}(k)+a_{dq}(k)\over d}.$$
\end{Prop}
{\it Proof}. We consider the formula given in the previous corollary. The divisors of
$M_k/2$ are either of the form $d$ with $d|kR_k$ or of the form $dk^2$ with $d|R_k$. For
the latter $d$ we find, using (\ref{naarkwadraatvrijekern}) that
$a_{dk^2}(k)=a_{dk}(1)=-\mu(dk)$ and hence
$\sum_{d|M_k/2}(a_d(k)+a_{dq}(k))/d=\sum_{d|kR_k}(a_d(k)+a_{dq}(k))/d$. Now
suppose that $d|R_k$. Using Proposition \ref{simpelzeg} we infer that
$$a_{dk}(k)+a_{dkq}(k)=a_d(k)+a_{dq}(k)+\mu(d)+\mu(dq)=a_d(k)+a_{dq}(k).$$
Using this observation it follows that
$$\sum_{d|kR_k}{a_d(k)+a_{dq}(k)\over d}=(1+{1\over k})\sum_{d|R_k}{a_d(k)+a_{dq}(k)\over d},$$
whence the result follows. \qed\\

\noindent Clearly $\delta(a_n(k)=0)=1-\sum_{v\ne 0}\delta(a_n(k)=v)$, where the latter
sum has only finitely many non-zero values and it is a finite computation to determine those $v$ for
which $a_n(k)=v$ for some $n$.   
For $1\le k\le 16$ the non-zero values of $\zeta(2)\delta(a_n(k)=v)$ are given
in Table 4 (except for $v=0$). A more extensive version of Table 4 is provided
by Table 11.\\

\centerline{{\bf Table 4:} Value  of $\zeta(2)\delta(a_n(k)=v)$}
\begin{center}
\begin{tabular}{|c|c|c|c|c|}
\hline
  & $v=-2$ & $v=-1$ & $v=1$ & $v=2$ \\
\hline
$k=1$  & $0$ & $1/2$ & $1/2$ & $0$ \\
\hline
$k=2$ & $0$ & $1/12$ & $7/12$ & $0$ \\
\hline
$k=3$ & $0$ & $5/24$ & $3/8$ & $0$ \\
\hline
$k=4$  & $0$ & $1/6$ & $1/2$ & $0$ \\ 
\hline
$k=5$  & $0$ & $13/80$ & $23/80$ & $0$ \\
\hline
$k=6$ & $0$ & $25/144$ & $67/144$ & $0$ \\
\hline
$k=7$ & $1/576$ & $577/2688$ & $731/2688$ & $1/1152$ \\
\hline
$k=8$  & $0$ & $1/8$ & $5/12$ & $0$ \\ 
\hline
$k=9$  & $0$ & $65/384$ & $347/1152$ & $0$ \\ 
\hline
$k=10$  & $0$ & $161/960$ & $347/960$ & $0$ \\ 
\hline
$k=11$  & $1/2304$ & $8299/50688$ & $11489/50688$ & $1/4608$ \\ 
\hline
$k=12$  & $0$ & $349/2304$ & $1009/2304$ & $0$ \\ 
\hline
$k=13$  & $43/48384$ & $219269/1257984$ & $277171/1257984$ & $43/96768$ \\ 
\hline
$k=14$  & $13/21504$ & $2395/21504$ & $2319/7168$ & $1/2304$ \\ 
\hline
$k=15$  & $13/32256$ & $1345/7168$ & $97247/322560$ & $13/64512$ \\ 
\hline
$k=16$  & $5/21504$ & $12149/64512$ & $1127/3072$ & $5/2688$ \\ 
\hline
\end{tabular}
\end{center}
\medskip

\noindent Let us now look at Theorem \ref{eenvoudig} and Theorem \ref{vier} from the viewpoint of
computational complexity. In Theorem \ref{eenvoudig} the sum has $p(k)$ terms and the famous
estimate of Hardy and Littlewood yields that $\log p(k)\sim \pi\sqrt{2k/3}$ as $k$ tends to
infinity. In Theorem \ref{vier} we sum over $t(k)$ terms where
$\log t(k)\sim \pi(k)\log 2 \sim k\log 2/\log k$. So Theorem \ref{eenvoudig} yields the
computational superior method. Theorem \ref{vier} is, however, much more easily implemented.
Using Proposition \ref{simpelzeg} in combination with Theorem \ref{vier} a result comparable
with Theorem \ref{eenvoudig} is obtained. Indeed, if one starts with Proposition \ref{priempiet}
and invokes Proposition \ref{simpelzeg}, one obtains a sum over partitions of $k$, where now
$2$ for example does not occur in the partition. This yields a result superior in complexity
to that provided by Theorem \ref{eenvoudig}.\\
\indent As already pointed out by M\"oller one could use his method to study the value distribution of
$a_n(k)$ in case e.g. $B(k)=2$ by considering the integer $a_n(k)(a_n(k)-1)/2$ to determine
$\delta(a_n(k)=-1)$ for example. This then yields a sum with $p(k)^2$ terms and this results
in an algorithm that has worse complexity than that provided by Theorem \ref{vier}. Aside from this,
this seems to be, from the practical point of view, an unwieldy method. 

\subsubsection{Some observations related to Table 4}
In this section we make some observations regarding Table 4 and Table 11, which is a much
more extensive version of Table 4.\\
\indent Let us put ${\cal B}(k)=\{a_n(k):n\in \mathbb N\}$. Recall that $|{\cal B}(k)|=B(k)$.
Recall also that $-1,0,1\in {\cal B}(k)$. Using this and Corollary \ref{corvier} we infer
that if $k$ is odd, then ${\cal B}(k)$ is symmetric: we have $v\in {\cal B}(k)$ iff
$-v\in {\cal B}(k)$. For $k$ is even numerical results suggest that often ${\cal B}(k)$ is
not symmetric. For $k\le 75$ it is true that if $v\in {\cal B}(k)$ and $v$ is negative, then
$-v\in {\cal B}(k)$. An other observation that can be made is that for $k\le 75$ it is
true that ${\cal B}(k)$ consists of consecutive integers, i.e. if $v_0<v_1$ are
in ${\cal B}(k)$, then so are all integers between $v_0$ and $v_1$.\\
\indent Let us define  ${\cal B}_j(k)=\{a_n(k):n\equiv j({\rm mod~}2)\}$, for $1\le j\le 2$.
It is easy to see that always $0\in {\cal B}_j(k)$.
Note that for $k$ is even we have ${\cal B}_0(k)\subseteq {\cal B}_1(k)={\cal B}(k)$. For
$k$ is odd we have $v\in {\cal B}_1(k)$ iff $-v \in {\cal B}_0(k)$. We could also
express this as ${\cal B}_0(k)=\{-v:v\in B_1(k)\}$. Thus in this case, if we know
${\cal B}_0(k)$, we also know ${\cal B}_1(k)$.
Inspection of Table 11 shows that for odd integers $k$ with $B(k)\ge 2$ often 
$\delta(a_n(k)=B(k))$ and $\delta(a_n(k)=-B(k))$ differ by
a factor two. Regarding this situation we have the following
result:
\begin{Prop}
\label{factortwo}
Let $k$ be odd and $v\ne 0$. We have 
$$2\delta(a_n(k)=v)=\delta(a_n(k)=-v){\rm ~iff~}v\not\in {\cal B}_1(k) 
~({\rm that~is~iff~} -v\not\in {\cal B}_0(k)).$$
\end{Prop}
{\it Proof}. `$\Rightarrow$' An easy consequence of Corollary \ref{corvier}. `$\Leftarrow$' This
uses in addition that $\{a_n(k):2\nmid n\}=\{a_d(k):d|M_k/2\}$. \qed\\

\noindent Example. Inspection of Table 4 shows that the condition of the proposition is 
satisfied for $k=7$ and $v=2$. It thus follows that there is no even integer $n$ for which
$a_n(7)=-2$ (whereas $a_{105}(7)=-2$). Further examples can be derived from Table 5.\\

\vfil\eject
\centerline{{\bf Table 5:} Set theoretic difference ${\cal B}(k)\backslash {\cal B}_0(k)$}
\medskip
\begin{center}
\begin{tabular}{|c|c|c|c|c|}
\hline
$k=7$    & $\{-2\}$ & $k=11$ & $\{-2\}$ \\
\hline
$k=13$ & $\{-2\}$ & $k=15$ & $\{-2\}$ \\
\hline
$k=17$ & $\{-3\}$ & $k=19$ & $\{-3\}$ \\
\hline
$k=21$  & $\{-3\}$ & $k=23$ & $\{-4,-3\}$ \\
\hline
$k=25$  & $\{-3\}$ & $k=31$ & $\{-4\}$ \\
\hline
$k=35$  & $\{5\}$ & $k=37$ & $\{5\}$ \\
\hline
$k=39$  & $\{5,6\}$ & $k=43$ & $\{-7\}$ \\
\hline
$k=45$  & $\{-7\}$ & $k=47$ & $\{-9,-8\}$ \\
\hline
$k=51$  & $\{8\}$ & $k=53$ & $\{9,10,11,12,13\}$ \\
\hline
\end{tabular}
\end{center}
\medskip

\noindent For Proposition \ref{factortwo} to be of some mathematical value we would hope that
  infinitely often ${\cal B}_0(k)$ is strictly contained in ${\cal B}(k)$. The next result
 shows that $a_{2n}(k)$ assumes all values as $n$ and $k$ range over the integers. 
\begin{Prop}
\label{suzukieven}
We have $\{a_{2n}(k):n,k\in \mathbb N\}=\mathbb Z$.
\end{Prop}
{\it Proof}. Given any integer $s\ge 2$ it is a consequence of the prime number theorem that
it is possible to find primes $2<p_1<p_1+2<p_2< \cdots < p_s \le p_1+p_2-2$. Let $k=p_1+p_2$
and let $q$ be any prime $>k$. If $s$ is even, let $m=2p_1\cdots p_sq$, otherwise
let $m=2p_1\cdots p_s$. We claim that
\begin{equation}
\label{suzieven}
\Phi_m(X)=\sum_{j=0}^{k-3}a_m(j)X^j+(s-1)X^{k-2}-(s-1)X^{k-1} ~({\rm mod~}X^k).
\end{equation}
The claim together with $a_4(1)=0$ yields the result. We next prove (\ref{suzieven}). Using
the observation that $\{d:d|m,~d<k\}=\{1,2,p_1,\cdots,p_s,2p_1\}$ we infer 
from (\ref{thanga}) that mod $X^k$ we
have
\begin{eqnarray}
\Phi_m(X)&\equiv &{(1-X)(1-X^{2p_1})\over (1-X^2)(1-X^{p_1})\cdots (1-X^{p_s})}~({\rm mod~}X^k)\cr
&\equiv & {1+X^{p_1}\over (1+X)(1-X^{p_2})\cdots (1-X^{p_s})}~({\rm mod~}X^k)\cr
&\equiv & (1-X+X^2-X^3\cdots)(1+X^{p_1}+\cdots+X^{p_s})~({\rm mod~}X^k)\cr
&\equiv &\sum_{j=0}^{k-3}a_m(j)X^j+(s-1)X^{k-2}-(s-1)X^{k-1} ~({\rm mod~}X^k).\nonumber
\end{eqnarray}
This concludes the proof. \qed\\

\noindent The above argument shows for example that 
$$\Phi_{2\cdot 19\cdot 23\cdot 29\cdot 31\cdot 37\cdot 39}(X)=\cdots-5X^{40}+5X^{41}+\cdots$$
\section{Value distribution of $s_m(p)$ and $S_m(p)$}
\subsection{Ramanujan sums, cyclotomic polynomials and primitive roots}
\noindent Proposition \ref{primitivetounity} connects the mod $p$ reductions of totally symmetric functions in primitive roots
to totally symmetric functions in primitive roots of unity. Our proof rests on the following simple result.
\begin{Prop}
\label{flauwflauw}
Let $p$ be a prime not dividing $m$. Then $p|\Phi_m(a)$ iff the order of $a({\rm mod~}p)$ is $m$.
\end{Prop}
{\it Proof}. If the order of $a({\rm mod~}p)$ equals $m$, then $a^d\not\equiv 1({\rm mod~}p)$ for every $d<m$
and $a^m\equiv 1({\rm mod~}p)$. By (\ref{basiccyclo}) we then infer that $p|\Phi_m(a)$. Since there are $\varphi(m)$ elements
of order $m$ in $(\mathbb Z/p\mathbb Z)^*$, the equation $\Phi_m(X)\equiv 0({\rm mod~}p)$ has at least
$\varphi(m)$ distinct solutions mod $p$, however since the degree of $\Phi_m(X)$ equals $\varphi(m)$, these
are all the solutions. \qed\\

\noindent Remark. We think this proof is more enlighting than that given in \cite[p. 209]{Murty}.
A related, less well-known result is that $a(a^{p-1}-1)\Psi'_{p-1}(a)/\Psi_{p-1}(a)$ is
congruent to $-1({\rm mod~}p)$ if $a$ is a primitive root mod $p$ and divisible by $p$ otherwise \cite{Nicol}.

\begin{Prop}
\label{primitivetounity}
Let $p$ be a prime. Put $t=\varphi(p-1)$. 
Let $g_1,\cdots,g_t$ be the modulu $p$ distinct primitive roots.
Let $1\le j_1<j_2<\cdots <j_t\le p-1$ be the natural numbers coprime to $p-1$.
Let $f(y_1,\cdots,y_t)$ be any totally symmetric function in the variables
$y_1,\cdots,y_t$.
Then 
$$f(g_1,\cdots,g_t)\equiv f(\zeta_{p-1}^{j_1},\cdots,\zeta_{p-1}^{j_t}) ~({\rm mod~}p).$$
\end{Prop}
{\it Proof}. By Newton's result on totally symmetric polynomial functions it is sufficient
to prove the result for all elementary totally symmetric polynomial functions.
By Proposition \ref{flauwflauw} we infer that
\begin{equation}
\label{factor1}
\Phi_{p-1}(x)\equiv (x-g_1)\cdots (x-g_t)~({\rm mod~}p).
\end{equation}
On the other hand we have
\begin{equation}
\label{factor2}
\Phi_{p-1}(x)=(x-\zeta_{p-1}^{j_1})\cdots (x-\zeta_{p-1}^{j_t}).
\end{equation}
From (\ref{factor1}) and (\ref{factor2}), the result easily follows. \qed\\

\noindent Let $g_1,\cdots,g_{\varphi(p-1)}$ denote the distinct primitive roots mod $p$ with 
$1\le g_j\le p-1$.
Recall that $$S_k(p)=\sum_{i=1}^{\varphi(p-1)}g_i^k,$$
and that $s_k(p)$ is the $k$th elementary totally symmetric function in $g_1,\cdots,g_{\varphi(p-1)}$.
By  Proposition \ref{primitivetounity} we infer that
$$S_k(p)\equiv c_{p-1}(k) ~({\rm mod~}p).$$
The latter congruence together with part 2 of Proposition \ref{basicramanujan} 
then yields the following lemma.
\begin{Lem} 
\label{moller}
Let $p$ be a prime. Then
$$S_k(p)\equiv \mu\left({p-1\over (p-1,k)}\right){\varphi(p-1)\over \varphi({p-1\over (p-1,k)})} ~({\rm  mod~}p).$$
\end{Lem}
This result is due to
Moller \cite{Moller}, who proved it in a different (and more roundabout) way. 
Note that in case $k=1$ we have $S_1(p)\equiv \mu(p-1)~{\rm (mod~}p)$, a 
result first established by Gauss \cite{Gauss}.\\

\noindent If $\varphi(p-1)<k$, then $s_k(p)=0$. If $k=\varphi(p-1)$, then we have, by Proposition 
\ref{primitivetounity},
$$s_k(p)=g_1g_2\cdots g_{\varphi(p-1)}\equiv (-1)^{\varphi(p-1)}\Phi_{p-1}(0)\equiv 1({\rm mod~}p).$$
If $\varphi(p-1)>k$, then $s_k(p)\equiv (-1)^ka_{p-1}(\varphi(p-1)-k)=(-1)^k a_{p-1}(k)({\rm mod~}p)$. Taking
the results of these three cases together we infer that
$$s_k(p)\equiv (-1)^ka_{p-1}(k) ({\rm mod~}p).$$ 
The latter congruence and the congruence $S_k(p)\equiv c_{p-1}(k)~({\rm mod~}p)$ are
the starting point of our analysis
of their value distribution which is carried out in the next two sections.

\subsection{Value distribution of $S_k(p)$ mod $p$}
As we have seen $S_k(p)\equiv c_{p-1}(k)~({\rm mod~}p)$ and so instead of the value distribution of
$S_k(p)$ mod $p$ we can investigate the value distribution of $c_{p-1}(k)$.
Given a natural number $k$, let $D_k=\{q_1,\dots,q_t\}$ denote the set of prime divisors of $k$, ordered in size. We have
$t=\# D_k=\omega(k)=\sum_{p|k}1$. The following observation will play an important role.
\begin{Prop}
\label{vastewaarde}
On each prime $p$ with $\nu_{D_k}(p-1)$ prescribed and $\mu_{D_k}(p-1)$ (respectively $\mu_{D_k}(p-1)^2$)  
prescribed, $c_{p-1}(k)$ (respectively $|c_{p-1}(k)|$) assumes the same value. If $\mu_{D_k}(p-1)=0$, this
value is zero.
\end{Prop}
Our most general result regarding the value distribution of $c_{p-1}(k)$ reads as follows.
\begin{Thm} 
\label{v}
Let $v\ne 0$.\\
{\rm 1)} For every $H>0$ we have
$$\sum_{p\le x \atop |c_{p-1}(k)|=v}1=c_v A{\rm Li}(x)+ O\left({x\over \log^H x}\right),$$
where the $O$-constant depends at most on $k$ and $H$ and $c_v$ is a rational number.
{\rm 2)} If we assume in addition Conjecture {\rm 1}, then the density of primes $p$ with $c_{p-1}(k)=v$ exists
and equals $c_vA/2$.
\end{Thm}
{\it Proof}. Part 1. Write $\nu_{D_k}(p-1)=(e_1,\dots,e_{\omega(k)})$. If $e_j\ge \nu_{q_j}(k)+2$ for some
$1\le j\le \omega(k)$, then $\mu((p-1)/(p-1,k))=0$ and hence $c_{p-1}(k)=0$. It follows there are only finitely
many possibilities for $\nu_{D_k}(p-1)$ to be considered. To each of them, by Proposition \ref{vastewaarde},  we
can associate an unique non-zero value of $c_{p-1}(k)$. We then have
$$\sum_{p\le x\atop |c_{p-1}(k)|=v}1=\sum_{p\le x\atop \nu_{D_k}(p-1)\in T}\mu_{D_k}(p-1)^2,$$
for some effectively computable set $T$. The result then follows by Proposition \ref{mirsky2}. 
Part 2. The proof of part 2 is similar to that of part 1 and left to the reader. \qed\\

\noindent Conjecture 1 together with Theorem \ref{mirsky}
(with $r=1$ and $k=2$) and the observation that 
$c_{p-1}(1)=\mu(p-1)$, then yields the following result demonstrating part 2 of Theorem \ref{v} (with $k=1$).
\begin{Prop} 
\label{stellingtwee}
Assuming Conjecture {\rm 1} we have
$$\delta(c_{p-1}(1)=j)=\cases{{A\over 2} &if $j=-1$;\cr
1-A &if $j=0$;\cr
{A\over 2} &if $j=1$.}
$$
\end{Prop}
\noindent Numerically we find up to the first $10^6$ primes `densities' of
$0.18732$, $0.625881$ and $0.186799$ respectively (for $-1,0,1$). This should be compared
to the conjectural values $0.1869\cdots$, $0.6260\cdots$ and $0.1869\cdots$, respectively.\\

\noindent We demonstrate part 1 of Theorem \ref{v} in Table 6 for $k=15$, where we took
$x$ such that $\pi(x)=10^6$ and rounded the result in the column
`numerical' to the sixth decimal.\\

\centerline{{\bf Table 6:} Density of values of $|S_{15}(p)|$ mod $p$}
\begin{center}
\begin{tabular}{|c|c|c|c|}
\hline
$v$  & {\rm theoretical} & {\rm numerical} & 
{\rm approximate}\\
\hline
$0$  & $1-561A/475$ & $0.558339$ & $0.558178$ \\
\hline
$1$  & $9A/19$ & $0.177137$ & $0.177157$ \\
\hline
$2$  & $6A/19$ & $0.118091$ & $0.118138$ \\
\hline
$3$  & $2A/19$ & $0.039364$ & $0.039353$ \\
\hline
$4$  & $12A/95$ & $0.047237$ & $0.047316$ \\
\hline
$5$  & $12A/475$ & $0.009447$ & $0.009457$ \\
\hline
$8$  & $8A/95$ & $0.031491$ & $0.031508$ \\
\hline
$10$  & $8A/475$ & $0.006298$ & $0.006317$ \\
\hline
$12$  & $8A/285$ & $0.010497$ & $0.010486$ \\
\hline
$15$  & $8A/1425$ & $0.002100$ & $0.002090$ \\
\hline
\end{tabular}
\end{center}
\medskip

\noindent Let $r$ be a prime and let $k=r^e$ for some $e\ge 0$. By Lemma \ref{moller} and (\ref{varquotient}) we deduce 
\begin{equation}
\label{kispriem}
c_{p-1}(k)=\cases{\mu(p-1) &if $\nu_r(p-1)=0$;\cr
r^f(1-{1\over r})\mu({p-1\over r^f}) &if $\nu_r(p-1)=f$, $1\le f\le e$;\cr
r^e\mu({p-1\over r^e}) &if $\nu_r(p-1)>e$.}
\end{equation}
As regards to the average size of $c_{p-1}(k)$ we have the following result.
\begin{Thm}
\label{absoluteaverage}
Let $k$ be any natural number, and
$H$ any positive real number. Then
$$\sum_{p\le x}|c_{p-1}(k)|=
A\prod_{q|k}\left(1+{\nu_q(k)(q-1)^2\over q^2-q-1}\right){\rm Li}(x)
+O\left({x\over \log^H x}\right),$$
where the $O$-constant depends at most on $k$ and $H$.
\end{Thm}
\begin{cor}
The average of $|c_{p-1}(k)|/A$ over the primes $p$ is a multiplicative function in $k$.
\end{cor}
Thus the average of $|c_{p-1}(k)|$ is a semi-multiplicative function in $k$. The analog of this
for the natural number case also holds true: $M_n(|c_n(m)|)$ is a semi-multiplicative function in
$m$ (we leave this as an exercise to the reader).\\
 
\noindent {\it Proof of Theorem} \ref{absoluteaverage}. Let us denote by $\delta(q,e)$ the density of primes $p$ such that
$\nu_{q}(p-1)=e$. By Lemma \ref{siegelwalfisz} one infers that
$$\delta(q,e)=\cases{1-{1\over q-1} &if $e=0$;\cr
q^{-e} &if $e\ge 1$.}$$ 
By property 6
of Proposition \ref{basicramanujan}, Proposition \ref{mirsky2}, the remark following 
Proposition \ref{mirsky2} and (\ref{kispriem}), we infer
that the theorem holds with a constant given by 
$$\delta(\nu_{S_k}(p-1)\ne 0)\prod_{q|k}\left(\delta(q,0)+\sum_{1\le f\le \nu_q(k)}\varphi(q^f)\delta(q,f)
+q^{\nu_q(k)}\delta(q,\nu_q(k)+1)\right).$$
On evaluating this constant, the result follows. \qed\\

\noindent Theorem \ref{absoluteaverage} is demonstrated in Table 7 (with again $x$ such that
$\pi(x)=10^6$).\\

\centerline{{\bf Table 7:} Average of $|c_{p-1}(k)|$}
\begin{center}
\begin{tabular}{|c|c|c|c|}
\hline
$k$  & {\rm theoretical} & {\rm numerical} & 
{\rm approximate}\\
\hline
$8$  & $4A$ & $1.494779$ & $1.495823$ \\
\hline
$21$  & $693A/205$ & $1.264572$ & $1.264153$ \\
\hline
$24$  & $36A/5$ & $2.689772$ & $2.692482$ \\
\hline
$27$  & $17A/5$ & $1.272214$ & $1.271450$ \\
\hline
$30$  & $126A/19$ & $2.479323$ & $2.479918$ \\
\hline
$36$  & $39A/5$ & $2.917172$ & $2.916855$ \\
\hline

\hline
\end{tabular}
\end{center}
\medskip

\noindent Our most general result in this section concerns the moments of $c_{p-1}(k)$.
\begin{Thm}
\label{moment}
Let $k$ be a natural number.\\ 
{\rm 1)}. Let $z\ne 1$ be a positive real number and $H$ be a  positive real number. Then
$$\sum_{p\le x}|c_{p-1}(k)|^{z}=$$
$$A\prod_{q|k}\left(1+{(q^{\nu_q(k)(z-1)}-1)(q-1)[(q-1)^{z}+q^{z-1}-1]\over (q^2-q-1)(q^{z-1}-1)}\right){\rm Li}(x)
+O\left({x\over \log^H x}\right),$$
where the $O$-constant depends at most on $k$ and $H$.\\
{\rm 2)}. If we assume Conjecture {\rm 1} and $j$ is an odd natural number, then we have
 $\sum_{p\le x}c_{p-1}(k)^j=o(\pi(x))$, where the $o$-constant depends at most
on $k$ and $j$.
\end{Thm}
{\it Proof}. Part 1. Proceeding as in the proof of Theorem \ref{absoluteaverage} we infer that the result holds with
constant
$$\delta(\nu_{S_k}(p-1)\ne 0)\prod_{q|k}\left(\delta(q,0)+\sum_{1\le f\le \nu_q(k)}\varphi(q^f)^{2j}\delta(q,f)
+q^{2j\nu_q(k)}\delta(q,\nu_q(k)+1)\right),$$
which is easily seen to equal the claimed constant.\\
Part 2. Similar to the proof of Theorem \ref{v}. \qed\\

\noindent Remark. If in part 1 we take the limit for $z$ tending to one, then we obtain 
Theorem \ref{absoluteaverage}.

\subsection{Value distribution of $s_k(p)$ mod $p$}
\label{vdofsk}
Recall that $s_k(p)\equiv (-1)^ka_{p-1}(k) ({\rm mod~}p)$. It follows that
$|\{s_k(p) {\rm ~mod~}p\}|\le  |\{ a_n(k):n\ge 1\}|$. A presumably difficult question
is to investigate under which conditions equality holds here. A related question is what the set
$\{a_{p-1}(k): k\ge 1,~p{\rm ~prime~}\}$ looks like (cf. Suzuki's result in Section 4.1).\\ 

\noindent First we consider the value distribution of $s_k(p)$ for
some small values of $k$. When $k=1$ we have $s_1(p)=S_1(p)$ and hence
for this case we refer to the previous section.
Note that $s_2(p)=\sum_{1\le i<j\le \varphi(p-1)}g_ig_j$. 
In the remainder of
this section we assume that $p>2$ (some of the assertions will also be valid for $p=2$).
Using Lemma \ref{moller} we
infer that
$$s_2(p)={1\over 2}(S_1(p)^2-S_2(p))\equiv {1\over 2}\left((\mu(p-1))^2-{\varphi(p-1)\over \varphi({p-1\over 2})}
\mu({p-1\over 2})\right) ({\rm mod~}p);$$
where the latter identity is Theorem 3 of \cite{DD}. By (\ref{varquotient}) we have
$$
{\varphi(p-1)\over \varphi({p-1\over 2})}=
\cases{
2 &if $p\equiv 1({\rm mod~}4)$;\cr
1 &if $p\equiv 3({\rm mod~}4)$.}
$$
On using this we infer that
\begin{equation}
\label{stweep}
s_2(p)\equiv
\cases{
-\mu({p-1\over 2})  &if $p\equiv 1({\rm mod~}4)$;\cr
\mu(p-1)(\mu(p-1)+1)/2 &if $p\equiv 3({\rm mod~}4)$.}
\end{equation}
(For notational convenience we will write $s_k(p)\equiv a$ for $s_k(p) \equiv a({\rm mod~}p)$.)
Since $s_2(2)=0$ it follows that the reduction of
$s_2(p)$ mod $p$ is in $\{-1,0,1\}$. The following assertion
conjecturally resolves the question of the Dence brothers
stated in the introduction.
\begin{Prop}
We have, assuming the validity of Conjecture {\rm 1},
$$
\delta(s_2(p)\equiv j)=
\cases{
{A\over 4} &if $j=-1$;\cr
1-A &if $j=0$;\cr
{3\over 4}A &if $j=1$.}
$$
\end{Prop}
{\it Proof}. Equation (\ref{stweep}) suggests to consider the 
primes $p$ with $p\equiv 1({\rm mod~}4)$ and $p\equiv 3({\rm mod~}4)$
separately. Note that for $x\ge 2$,
$$\sum_{p\le x,~p\equiv 3({\rm mod~}4)\atop \mu(p-1)=1}1=-1+\sum_{p\le x\atop \mu(p-1)=1}1.$$
The latter quantity we (conditionally) evaluated in Proposition \ref{stellingtwee}. It follows that the density of primes $p$ such that
$p\equiv 3({\rm mod~}4)$ and $s_2(p) \equiv 1$ equals $A/2$ and that the density of primes $p$ such that
$p\equiv 3({\rm mod~}4)$ and $s_2(p) \equiv 0$ equals $1/2-A/2$.\\
\indent The case $p\equiv 1({\rm mod~}4)$ can be dealt with Proposition 
\ref{mirsky2} (with $t=1$, $q_1=2$ and $e_1=2$) and Conjecture 1. The results are summed up in 
Table 8.\qed\\

\centerline{{\bf Table 8:} Value distribution of $s_2(p)$ mod $p$}
\begin{center}
\begin{tabular}{|c||c|c|c|c|}
\hline
$\delta$  & $-1$ & $0$ & 
$1$ & $+$\\
\hline
$\nu_2(p-1)\le 1$  & $0$ & $1/2-A/2$ & $A/2$ & $1/2$ \\
{\rm numerical}  & $0.000000$ & $0.313022$ & $0.186978$ & $0.500000$ \\
$\pi(x)=10^6$    & $0.000000$ & $0.313403$ & $0.186798$ & $0.500201$ \\
\hline
$\nu_2(p-1)=2$ & $A/4$ & $1/2-A/2$ & $A/4$ & $1/2$ \\
{\rm numerical}  & $0.093489$ & $0.313022$ & $0.093489$ & $0.500000$ \\
$\pi(x)=10^6$ & $0.093939$ & $0.312813$ & $0.093047$ & $0.499799$ \\
\hline
$+$ & $A/4$ & $1-A$ & $3A/4$ & $1$ \\
{\rm numerical}  & $0.093489$ & $0.626044$ & $0.280467$ & $1.000000$ \\
 $\pi(x)=10^6$  &  $0.093939$ & $0.626216$ & $0.279845$ & $1.000000$ \\

\hline
\end{tabular}
\end{center}
\medskip

\noindent By Newton's formula $s_k(p)$ can be expressed as a polynomial
in $S_r(p)$ with $1\le r\le k$. To be precise, we have
\begin{equation}
\label{newton}
(-1)^ks_k(p)=\sum{(-1)^{k_1+k_2+\cdots}\over k_1!1^{k_1}k_2!2^{k_2}\cdots}S_1(p)^{k_1}S_2(p)^{k_2}\cdots,
\end{equation}
where the sum extends over all solutions $(k_1,k_2,\cdots)$ of $k_1+2k_2+\cdots=k$. 
Next we consider $s_3(p)$. Taking $k=3$ in (\ref{newton}) we compute
$$s_3(p)={S_1(p)^3+2S_3(p)-3S_1(p)S_2(p)\over 6}.$$ 
Invoking Lemma \ref{moller} a more explicit formula for
$s_3(p)$ can then be derived. With $\beta=\nu_3(p-1)$ we find
$$
s_3(p)=
\cases{
\mu(p-1)(\mu(p-1)+1)/2  &if $\beta=0$;\cr
\mu(p-1)(\mu(p-1)-1)/2 &if $\beta=1$;\cr
\mu({p-1\over 3}) &if $\beta\ge 2$.}
$$
In each of these cases Proposition \ref{mirsky2} and Conjecture 1 yield the densities with which the values
$-1,0$ and $1$ are assumed. The results are summed up in Table 9.\\
\vfil\eject
\centerline{{\bf Table 9:} Value distribution of $s_3(p)$}
\begin{center}
\begin{tabular}{|c||c|c|c|c|}
\hline
$\delta$  & $-1$ & $0$ & 
$1$ & $+$\\
\hline
$\beta=0$  & $0$ & $1/2-3A/10$ & $3A/10$ & $1/2$ \\
{\rm numerical}  & $0.000000$ & $0.387813$ & $0.112187$ & $0.500000$ \\
$\pi(x)=10^6$    & $0.000000$ & $0.388136$ & $0.112035$ & $0.500171$ \\
\hline
$\beta=1$ & $0$ & $1/3-A/5$ & $A/5$ & $1/3$ \\
{\rm numerical}  & $0.000000$ & $0.258542$ & $0.074791$ & $0.333333$ \\
$\pi(x)=10^6$ & $0.000000$ & $0.258378$ & $0.074883$ & $0.333261$ \\
\hline
$\beta\ge 2$  & $A/15$ & $1/6-2A/15$ & $A/15$ & $1/6$ \\
{\rm numerical}  & $0.024930$ & $0.116806$ & $0.024930$ & $0.166666$ \\
$\pi(x)=10^6$ & $0.025031$ & $0.116729$ & $0.024809$ & $0.166569$ \\
\hline
$+$ & $A/15$ & $1-19A/30$ & $17A/30$ & $1$ \\
{\rm numerical}  & $0.024930$ & $0.763161$ & $0.211908$ & $1.000000$ \\
 $\pi(x)=10^6$  & $0.025030$ & $0.763243$ & $0.211727$ & $1.000000$ \\

\hline
\end{tabular}
\end{center}
\medskip

\noindent Adding up the various contributions we infer the following result.
\begin{Prop}
We have, assuming the validity of Conjecture {\rm 1},
$$
\delta(s_3(p)\equiv j)=
\cases{
{A\over 15} &if $j=-1$;\cr
1-{19\over 30}A &if $j=0$;\cr
{17\over 30}A &if $j=1$.}
$$
\end{Prop}

\noindent Similarly, for $s_4(p)$ we find, with $\alpha=\nu_2(p-1)$ and $\beta=\nu_3(p-1)$,
$$
s_4(p)\equiv 
\cases{
\mu(p-1)(\mu(p-1)+1)/2  &if $\alpha=1$ and $\beta=0$;\cr
\mu(p-1)(1-\mu(p-1))/2 &if $\alpha=1$ and $\beta\ge 1$;\cr
\mu({p-1\over 2})(\mu({p-1\over 2})+1)/2  &if $\alpha=2$;\cr
-\mu({p-1\over 4}) &if $\alpha\ge 3$.}
$$

\noindent Proceeding as before we infer the following result.
\begin{Prop}
We have, assuming the validity of Conjecture {\rm 1},
$$
\delta(s_4(p)\equiv j)=
\cases{
{13\over 40}A &if $j=-1$;\cr
1-A &if $j=0$;\cr
{27\over 40}A &if $j=1$.}
$$
\end{Prop}

\noindent The above approach of dealing with the 
value distribution of $s_k(p)({\rm mod~}p)$ for small
$p$ clearly becomes more laborious as $k$ increases. The following
result is more systematic; it is analogous to Theorem \ref{vier}.
\begin{Thm}
Let $k\ge 1$ be fixed. 
Let $M_k=k\prod_{p\le k}p$ and let
$r>k$ be any prime. Then
$$\sum_{p\le x}a_{p-1}(k)={Ax\over 2\log x}
\prod_{2<q\le k}{q(q-2)\over q^2-q-1}\sum_{d|M_k\atop 2|d}{a_{d}(k)+a_{dr}(k)\over d}\prod_{q|d\atop q>2}
{q-1\over q-2}+o({x\over \log x}).$$
If $v\ne 0$, then
$$\sum_{p\le x\atop a_{p-1}(k)=v}1=B_2(k){x\over \log x}+o({x\over \log x}),$$
where  
$$B_2(k)={A\over 2}\prod_{2<q\le k}{q(q-2)\over q^2-q-1}
\left(\sum_{d|M_k,~2|d\atop a_{d}(k)=v}{1\over d}\prod_{q|d\atop q>2}
{q-1\over q-2}+\sum_{d|M_k,~2|d\atop a_{dr}(k)=v}{1\over d}\prod_{q|d\atop q>2}
{q-1\over q-2}\right).$$
In both cases the $o$-constant depends at most on $k$.
\end{Thm}
{\it Proof}. 
Put $N_k={\rm lcm}(1,2,\cdots,k)\prod_{p\le k}p$.
Let $S(k)=\{q_1,\dots,q_{\pi(k)}\}$ be the set of primes not exceeding $k$ and
$r_2=N_k\prod_{p\le k}p$. Furthermore, let $d|N_k$. We first
consider $\sum_{p\le x,~(p-1,r_2)=d}\mu((p-1)/d)$. Let $q\in S(k)$ and suppose that $p$ satisfies
$(p-1,r_2)=d$. Since $\nu_q(r_2)>\nu_q(d)$, it follows that $\nu_q(p-1)=\nu_q(d)$. Note that
$\mu((p-1)/d)=\mu_{S(k)}(p-1)$. Conjecture 1 implies now that 
\begin{equation}
\label{eerstem}
\sum_{p\le x,~(p-1,r_2)=d}\mu({p-1\over d})=\sum_{p\le x\atop \nu_{S(k)}(p-1)=(\nu_{q_1}(d),\dots,
\nu_{q_{\pi(k)}}(d))}\mu_{S(k)}(p-1)=o(\pi(x)),
\end{equation} 
where the o-constant depends at most on $k$.\\
\indent Next we consider
\begin{equation}
\label{tweedem}
\sum_{p\le x\atop (p-1,r_2)=d}\mu({p-1\over d})^2=\sum_{p\le x\atop \nu_{S(k)}(p-1)=(\nu_{q_1}(d),\dots,
\nu_{q_{\pi(k)}}(d))}\mu_{S(k)}(p-1)^2.
\end{equation}
By Proposition \ref{mirsky2} we then find that
$$\sum_{p\le x\atop (p-1,r)=d}\mu({p-1\over d})^2={1\over d}\prod_{q\le k\atop q\nmid d}\left(1-{1\over q-1}\right)
\prod_{q>k}\left(1-{1\over q(q-1)}\right){x\over \log x}+o({x\over \log x}).$$
Using Proposition \ref{naareindig} we have
$$\sum_{p\le x}a_{p-1}(k)=\sum_{d|N_k}\sum_{p\le x\atop (p-1,r_2)=d}
\left(A_1(d)\mu({p-1\over d})^2+B_1(d)\mu({p-1\over d})\right).$$
From the latter formula in combination with (\ref{eerstem}) and (\ref{tweedem}), we infer that
$$\sum_{p\le x}a_{p-1}(k)={Ax\over \log x}\prod_{2<q\le k}{q(q-2)\over q^2-q-1}\sum_{d|N_k\atop 2|d}{A_1(d)\over d}\prod_{q|d\atop q>2}
{q-1\over q-2}+o({x\over \log x}).$$
On noting that $A_1(d)=(a_{d}(k)+a_{dr}(k))/2$ and
using (\ref{naarkwadraatvrijekern}), the proof of the first assertion is completed. 
The proof of the second assertion is a variation of the argument above and left
to the reader. \qed\\

\noindent The latter theorem is demonstrated in Table 10. Using that the sum of the densities equals
one, one can easily find the density of primes for which $a_{p-1}(k)=0$ from the table, so we
omitted this information.  Notice that
the averages given for $1\le k\le 4$ are the same as can be inferred from the earlier 
presented results on the value
distribution of $S_1(p),s_2(p),s_3(p)$ and $s_4(p)$ respectively. The averages given also seem
to be consistent with numerical simulations.\\

\centerline{{\bf Table 10:} Conjectural value  distribution of $\delta(a_{p-1}(k)=v)/A$}
\begin{center}
\begin{tabular}{|c|c|c|c|c|c|}
\hline
$j$  & $-2$ & $-1$ & $1$ & $2$ & {\rm Average}\\
\hline
$k=1$  & $0$ & $1/2$ & $1/2$ & $0$ & $0$ \\
\hline
$k=2$ & $0$ & $1/4$ & $3/4$ & $0$ & $1/2$ \\
\hline
$k=3$ & $0$ & $17/30$ & $1/15$ & $0$ & $-1/2$ \\
\hline
$k=4$  & $0$ & $13/40$ & $27/40$ & $0$ & $7/20$ \\ 
\hline
$k=5$  & $0$ & $69/90$ & $6/95$ & $0$ & $-3/10$ \\
\hline
$k=6$ & $0$ & $443/1140$ & $47/95$ & $0$ & $121/1140$ \\
\hline
$k=7$ & $0$ & $13989/54530$ & $358/1435$ & $24/3895$ & $1/190$ \\
\hline
$k=8$  & $0$ & $16703/62320$ & $35873/62320$ & $0$ & $1917/6232$\\ 
\hline
$k=9$  & $0$ & $31477/70110$ & $2129/35055$ & $0$ & $-9073/23370$ \\ 
\hline
$k=10$  & $0$ & $267/820$ & $505/1558$ & $0$ & $-23/15580$ \\ 
\hline
\end{tabular}
\end{center}

\noindent It is not difficult to see that $\delta(a_{p-1}(k)=v)>0$ if and only if 
$a_{p-1}(k)=v$ for some prime $p$. We infer from the table that there is no prime $p$
such that $a_{p-1}(7)=-2$, whereas by Table 4 it follows that $\delta(a_{n}(7)=-2)=0.001055\cdots$. 
Indeed, the example
following Proposition \ref{verdubbeling} shows that there are no even integers $n$ for
which $a_n(7)=-2$.

\section{Open problems}
During research usually more questions are being raised than being resolved and the present work
is no exception. We conclude by formulating some questions. Gallot's data show that
2, 5 and 9 hold true for $k\le 100$.\\
{\rm 1)} Is Conjecture 1 true ?\\
{\rm 2)} Is M\"oller's conjecture that $0\le e_k\le 1/2$ true ?\\
{\rm 3)} What is the behaviour of $e_k$ for large $k$ ?\\
{\rm 4)} Compute $M(e_k)$.\\
{\rm 5)} Is $e_kk\prod_{p\le k}(p+1)\in \mathbb Z$ ?\\
{\rm 6)} Determine $\{a_{p-1}(k):p{\rm ~prime},k\in \mathbb N\}$.\\
{\rm 7)} Does $a_n(k)=-v$, $v\le 0$, imply there is an integer $m$ with $a_m(k)=v$ ?\\
{\rm 8)} Is $|{\cal B}(k)\backslash {\cal B}_0(k)|\ge 1$ infinitely often ?\\
{\rm 9)} Is $\delta(a_n(k)=1)\ge \delta(a_n(k)=-1)$ ?\\
{\rm 10)} What is the behaviour of $\delta(a_n(k)\not\in \{-1,0,1\})$ for large $k$ ? 
\medskip

\noindent {\tt Acknowledgement}. We are very grateful to dr. Gallot for assistance in computing
some of the tables. We thank
Prof. T. Dence for sending \cite{DD} and Prof. Thangadurai for sending us his beautiful survey paper
\cite{Thanga}, Prof. Bachman for some helpful e-mail correspondence and dr. Tegelaar for making
some plots using Mathematica.


\vfil\eject

\centerline{{\bf Table 11:} Value distribution of $a_n(k)$ for $1\le k\le 61$}
\medskip
\noindent (The number between brackets is $k\prod_{p\le k}(p+1)e_k$.)
\medskip

\noindent $e_{1} = 0 (0.000000) [0]\\
V_1[-1] = 1/2 (0.500000)\\
V_1[+1] = 1/2 (0.500000)\\
\\
e_2 = 1/2 (0.500000) [3]\\
V_2[-1] = 1/12 (0.083333)\\
V_2[+1] = 7/12 (0.583333)\\
\\
e_3 = 1/6 (0.166667) [6]\\
V_3[-1] = 5/24 (0.208333)\\
V_3[+1] = 3/8 (0.375000)\\
\\
e_4 = 1/3 (0.333333) [16]\\
V_4[-1] = 1/6 (0.166667)\\
V_4[+1] = 1/2 (0.500000)\\
\\
e_5 = 1/8 (0.125000) [45]\\
V_5[-1] = 13/80 (0.162500)\\
V_5[+1] = 23/80 (0.287500)\\
\\
e_6 = 7/24 (0.291667) [126]\\
V_6[-1] = 25/144 (0.173611)\\
V_6[+1] = 67/144 (0.465278)\\
\\
e_7 = 1/18 (0.055556) [224]\\
V_7[-2] = 1/576 (0.001736)\\
V_7[-1] = 577/2688 (0.214658)\\
V_7[+1] = 731/2688 (0.271949)\\
V_7[+2] = 1/1152 (0.000868)\\
\\
e_8 = 7/24 (0.291667) [1344]\\
V_8[-1] = 1/8 (0.125000)\\
V_8[+1] = 5/12 (0.416667)\\
\\
e_9 = 19/144 (0.131944) [684]\\
V_9[-1] = 65/384 (0.169271)\\
V_9[+1] = 347/1152 (0.301215)\\
\\
e_{10} = 31/160 (0.193750) [1116]\\
V_{10}[-1] = 161/960 (0.167708)\\
V_{10}[+1] = 347/960 (0.361458)\\$
\\
\vfil\eject
\noindent $e_{11} = 1/16 (0.062500) [4752]\\
V_{11}[-2] = 1/2304 (0.000434)\\
V_{11}[-1] = 8299/50688 (0.163727)\\
V_{11}[+1] = 11489/50688 (0.226661)\\
V_{11}[+2] = 1/4608 (0.000217)\\
\\
e_{12} = 55/192 (0.286458) [23760]\\
V_{12}[-1] = 349/2304 (0.151476)\\
V_{12}[+1] = 1009/2304 (0.437934)\\
\\
e_{13} = 13/288 (0.045139) [56784]\\
V_{13}[-2] = 43/48384 (0.000889)\\
V_{13}[-1] = 219269/1257984 (0.174302)\\
V_{13}[+1] = 277171/1257984 (0.220330)\\
V_{13}[+2] = 43/96768 (0.000444)\\
\\
e_{14} = 61/288 (0.211806) [286944]\\
V_{14}[-2] = 13/21504 (0.000605)\\
V_{14}[-1] = 2395/21504 (0.111375)\\
V_{14}[+1] = 2319/7168 (0.323521)\\
V_{14}[+2] = 1/2304 (0.000434)\\
\\
e_{15} = 2287/20160 (0.113442) [164664]\\
V_{15}[-2] = 13/32256 (0.000403)\\
V_{15}[-1] = 1345/7168 (0.187640)\\
V_{15}[+1] = 97247/322560 (0.301485)\\
V_{15}[+2] = 13/64512 (0.000202)\\
\\
e_{16} = 733/4032 (0.181796) [281472]\\
V_{16}[-2] = 5/21504 (0.000233)\\
V_{16}[-1] = 12149/64512 (0.188322)\\
V_{16}[+1] = 1127/3072 (0.366862)\\
V_{16}[+2] = 5/2688 (0.001860)\\
\\
e_{17} = 667/8064 (0.082713) [2449224]\\
V_{17}[-3] = 5/580608 (0.000009)\\
V_{17}[-2] = 281/193536 (0.001452)\\
V_{17}[-1] = 2353487/19740672 (0.119220)\\
V_{17}[+1] = 1981753/9870336 (0.200779)\\
V_{17}[+2] = 197/96768 (0.002036)\\
V_{17}[+3] = 5/1161216 (0.000004)\\$
\\
\vfil\eject
\noindent $e_{18} = 79/336 (0.235119) [7371648]\\
V_{18}[-2] = 961/1161216 (0.000828)\\
V_{18}[-1] = 5575/43008 (0.129627)\\
V_{18}[+1] = 212369/580608 (0.365770)\\
V_{18}[+2] = 341/1161216 (0.000294)\\
V_{18}[+3] = 17/1161216 (0.000015)\\
\\
e_{19} = 55/1344 (0.040923) [27086400]\\
V_{19}[-3] = 19/2322432 (0.000008)\\
V_{19}[-2] = 67813/69672960 (0.000973)\\
V_{19}[-1] = 42731243/264757248 (0.161398)\\
V_{19}[+1] = 67086449/330946560 (0.202711)\\
V_{19}[+2] = 54641/69672960 (0.000784)\\
V_{19}[+3] = 19/4644864 (0.000004)\\
\\
e_{20} = 221/960 (0.230208) [160392960]\\
V_{20}[-3] = 1/221184 (0.000005)\\
V_{20}[-2] = 307/165888 (0.001851)\\
V_{20}[-1] = 1417037/11612160 (0.122030)\\
V_{20}[+1] = 8240789/23224320 (0.354834)\\
V_{20}[+2] = 361/645120 (0.000560)\\
\\
e_{21} = 8207/120960 (0.067849) [49635936]\\
V_{21}[-3] = 1/331776 (0.000003)\\
V_{21}[-2] = 138259/69672960 (0.001984)\\
V_{21}[-1] = 14659501/69672960 (0.210404)\\
V_{21}[+1] = 69503/248832 (0.279317)\\
V_{21}[+2] = 101363/69672960 (0.001455)\\
V_{21}[+3] = 1/663552 (0.000002)\\
\\
e_{22} = 8467/95040 (0.089089) [68277888]\\
V_{22}[-3] = 1/221184 (0.000005)\\
V_{22}[-2] = 7/18432 (0.000380)\\
V_{22}[-1] = 8101001/42577920 (0.190263)\\
V_{22}[+1] = 70728809/255467520 (0.276860)\\
V_{22}[+2] = 3751/2322432 (0.001615)\\
V_{22}[+3] = 19/1658880 (0.000011)\\
\\
e_{23} = 629/11520 (0.054601) [1049956992]\\
V_{23}[-4] = 1/7962624 (0.000000)\\
V_{23}[-3] = 107/7962624 (0.000013)\\
V_{23}[-2] = 578371/418037760 (0.001384)\\
V_{23}[-1] = 5016105893/38459473920 (0.130426)\\
V_{23}[+1] = 1789013287/9614868480 (0.186067)\\
V_{23}[+2] = 730129/836075520 (0.000873)\\
V_{23}[+3] = 107/15925248 (0.000007)\\
V_{23}[+4] = 1/15925248 (0.000000)\\
\\
e_{24} = 7327/24192 (0.302869) [6077306880\\
V_{24}[-3] = 853/69672960 (0.000012)\\
V_{24}[-2] = 247859/278691840 (0.000889)\\
V_{24}[-1] = 898117/9289728 (0.096679)\\
V_{24}[+1] = 55697567/139345920 (0.399707)\\
V_{24}[+2] = 8227/10321920 (0.000797)\\
V_{24}[+3] = 713/34836480 (0.000020)\\
\\
e_{25} = 1087/18144 (0.059910) [1252224000]\\
V_{25}[-3] = 25/7962624 (0.000003)\\
V_{25}[-2] = 81607/238878720 (0.000342)\\
V_{25}[-1] = 1375923029/8360755200 (0.164569)\\
V_{25}[+1] = 469535011/2090188800 (0.224638)\\
V_{25}[+2] = 88489/334430208 (0.000265)\\
V_{25}[+3] = 25/15925248 (0.000002)\\
\\
e_{26} = 234433/1572480 (0.149085) [3240801792]\\
V_{26}[-3] = 25/5308416 (0.000005)\\
V_{26}[-2] = 11149/15482880 (0.000720)\\
V_{26}[-1] = 837724879/7245987840 (0.115612)\\
V_{26}[+1] = 127807265/483065856 (0.264575)\\
V_{26}[+2] = 54667/69672960 (0.000785)\\
V_{26}[+3] = 1289/557383680 (0.000002)\\
\\
e_{27} = 33491/362880 (0.092292) [2083408128]\\
V_{27}[-3] = 811/185794560 (0.000004)\\
V_{27}[-2] = 311893/836075520 (0.000373)\\
V_{27}[-1] = 244503443/1672151040 (0.146221)\\
V_{27}[+1] = 399281329/1672151040 (0.238783)\\
V_{27}[+2] = 201713/836075520 (0.000241)\\
V_{27}[+3] = 139/61931520 (0.000002)\\
\\
e_{28} = 84047/483840 (0.173708) [4066530048]\\
V_{28}[-3] = 1073/69672960 (0.000015)\\
V_{28}[-2] = 202387/278691840 (0.000726)\\
V_{28}[-1] = 5770579/39813120 (0.144942)\\
V_{28}[+1] = 88832659/278691840 (0.318749)\\
V_{28}[+2] = 96629/139345920 (0.000693)\\
V_{28}[+3] = 1/241920 (0.000004)\\
V_{28}[+4] = 1/7741440 (0.000000)\\
\\
e_{29} = 5021/103680 (0.048428) [35225729280]\\
V_{29}[-4] = 1837/25082265600 (0.000000)\\
V_{29}[-3] = 215897/16721510400 (0.000013)\\
V_{29}[-2] = 56528869/50164531200 (0.001127)\\
V_{29}[-1] = 5288720741/41564897280 (0.127240)\\
V_{29}[+1] = 256004923169/1454771404800 (0.175976)\\
V_{29}[+2] = 48779039/50164531200 (0.000972)\\
V_{29}[+3] = 24511/1857945600 (0.000013)\\
V_{29}[+4] = 3883/50164531200 (0.000000)\\
\\
e_{30} = 45893/241920 (0.189703) [142745587200]\\
V_{30}[-3] = 55183/8360755200 (0.000007)\\
V_{30}[-2] = 2194267/2090188800 (0.001050)\\
V_{30}[-1] = 541742161/3344302080 (0.161990)\\
V_{30}[+1] = 5867204267/16721510400 (0.350878)\\
V_{30}[+2] = 34039/23887872 (0.001425)\\
V_{30}[+3] = 450059/16721510400 (0.000027)\\
V_{30}[+4] = 8107/8360755200 (0.000001)\\
V_{30}[+5] = 29/5573836800 (0.000000)\\
\\
e_{31} = 155/5376 (0.028832) [717382656000]\\
V_{31}[-4] = 341/7644119040 (0.000000)\\
V_{31}[-3] = 25652623/802632499200 (0.000032)\\
V_{31}[-2] = 1728495157/1605264998400 (0.001077)\\
V_{31}[-1] = 4263256278479/24881607475200 (0.171342)\\
V_{31}[+1] = 4956093793291/24881607475200 (0.199187)\\
V_{31}[+2] = 2463624983/1605264998400 (0.001535)\\
V_{31}[+3] = 44559467/802632499200 (0.000056)\\
V_{31}[+4] = 341/15288238080 (0.000000)\\
\\
e_{32} = 4381/17920 (0.244475) [6279166033920]\\
V_{32}[-4] = 143/5096079360 (0.000000)\\
V_{32}[-3] = 143131/19818086400 (0.000007)\\
V_{32}[-2] = 32276503/38220595200 (0.000844)\\
V_{32}[-1] = 45508919251/535088332800 (0.085049)\\
V_{32}[+1] = 24918610681/76441190400 (0.325984)\\
V_{32}[+2] = 667940893/267544166400 (0.002497)\\
V_{32}[+3] = 8249729/107017666560 (0.000077)\\
V_{32}[+4] = 1207481/178362777600 (0.000007)\\
\\
e_{33} = 294509/3193344 (0.092226) [2442775449600]\\
V_{33}[-4] = 20065/10701766656 (0.000002)\\
V_{33}[-3] = 80453581/802632499200 (0.000100)\\
V_{33}[-2] = 412044869/401316249600 (0.001027)\\
V_{33}[-1] = 2239952947681/17657914982400 (0.126853)\\
V_{33}[+1] = 240607474937/1103619686400 (0.218017)\\
V_{33}[+2] = 1136159459/802632499200 (0.001416)\\
V_{33}[+3] = 309052627/1605264998400 (0.000193)\\
V_{33}[+4] = 398297/107017666560 (0.000004)\\$
\\
\vfil\eject
\noindent $e_{34} = 268801/3525120 (0.076253) [2080906813440]\\
V_{34}[-4] = 143/5096079360 (0.000000)\\
V_{34}[-3] = 20977/5096079360 (0.000004)\\
V_{34}[-2] = 1313801/1189085184 (0.001105)\\
V_{34}[-1] = 1500096902327/9096501657600 (0.164909)\\
V_{34}[+1] = 135360848387/568531353600 (0.238089)\\
V_{34}[+2] = 80945987/33443020800 (0.002420)\\
V_{34}[+3] = 3130531/21403533312 (0.000146)\\
V_{34}[+4] = 1071667/267544166400 (0.000004)\\
V_{34}[+5] = 719/19818086400 (0.000000)\\
\\
e_{35} = 69809/870912 (0.080156) [2251759104000]\\
V_{35}[-5] = 323/2090188800 (0.000000)\\
V_{35}[-4] = 242639/25082265600 (0.000010)\\
V_{35}[-3] = 280214099/1605264998400 (0.000175)\\
V_{35}[-2] = 428425901/229323571200 (0.001868)\\
V_{35}[-1] = 15515446957/107017666560 (0.144980)\\
V_{35}[+1] = 120969894043/535088332800 (0.226075)\\
V_{35}[+2] = 2211435001/1605264998400 (0.001378)\\
V_{35}[+3] = 292252153/1605264998400 (0.000182)\\
V_{35}[+4] = 4691237/321052999680 (0.000015)\\
V_{35}[+5] = 323/1045094400 (0.000000)\\
\\
e_{36} = 41333/207360 (0.199330) [5759584911360]\\
V_{36}[-3] = 5072809/76441190400 (0.000066)\\
V_{36}[-2] = 31579379/26754416640 (0.001180)\\
V_{36}[-1] = 92196104039/535088332800 (0.172301)\\
V_{36}[+1] = 12419782037/33443020800 (0.371371)\\
V_{36}[+2] = 124255577/107017666560 (0.001161)\\
V_{36}[+3] = 78450899/535088332800 (0.000147)\\
V_{36}[+4] = 22711/1651507200 (0.000014)\\
V_{36}[+5] = 56381/178362777600 (0.000000)\\
V_{36}[+6] = 14633/535088332800 (0.000000)\\
\\
e_{37} = 19073/544320 (0.035040) [39542741729280]\\
V_{37}[-5] = 30631/3210529996800 (0.000000)\\
V_{37}[-4] = 56904067/61000069939200 (0.000001)\\
V_{37}[-3] = 60188399/1742859141120 (0.000035)\\
V_{37}[-2] = 76311280099/61000069939200 (0.001251)\\
V_{37}[-1] = 323928333998587/2257002587750400 (0.143521)\\
V_{37}[+1] = 80619151138297/451400517550080 (0.178598)\\
V_{37}[+2] = 73686178829/61000069939200 (0.001208)\\
V_{37}[+3] = 763875323/15250017484800 (0.000050)\\
V_{37}[+4] = 51705001/30500034969600 (0.000002)\\
V_{37}[+5] = 30631/1605264998400 (0.000000)\\$
\\
\vfil\eject
\noindent $e_{38} = 246539/1451520 (0.169849) [196855040655360]\\
V_{38}[-4] = 7955363/2541669580800 (0.000003)\\
V_{38}[-3] = 18756869/150617456640 (0.000125)\\
V_{38}[-2] = 1660704037/1694446387200 (0.000980)\\
V_{38}[-1] = 137023551787/1694446387200 (0.080866)\\
V_{38}[+1] = 1265500142927/5083339161600 (0.248951)\\
V_{38}[+2] = 4927088377/2541669580800 (0.001939)\\
V_{38}[+3] = 96190153/1270834790400 (0.000076)\\
V_{38}[+4] = 17066753/10166678323200 (0.000002)\\
V_{38}[+5] = 629/356725555200 (0.000000)\\
\\
e_{39} = 961879/12804480 (0.075121) [89355942789120]\\
V_{39}[-6] = 1/179159040 (0.000000)\\
V_{39}[-5] = 179/627056640 (0.000000)\\
V_{39}[-4] = 77322241/12200013987840 (0.000006)\\
V_{39}[-3] = 6774345563/61000069939200 (0.000111)\\
V_{39}[-2] = 2729296757/2033335664640 (0.001342)\\
V_{39}[-1] = 121853678359577/793000909209600 (0.153661)\\
V_{39}[+1] = 181373355248401/793000909209600 (0.228718)\\
V_{39}[+2] = 13458490453/10166678323200 (0.001324)\\
V_{39}[+3] = 1683824087/12200013987840 (0.000138)\\
V_{39}[+4] = 7100759/642105999360 (0.000011)\\
V_{39}[+5] = 179/313528320 (0.000001)\\
V_{39}[+6] = 1/89579520 (0.000000)\\
\\
e_{40} = 3360923/18385920 (0.182799) [223014717849600]\\
V_{40}[-5] = 198679/26754416640 (0.000007)\\
V_{40}[-4] = 2961557/4066671329280 (0.000001)\\
V_{40}[-3] = 5036034173/20333356646400 (0.000248)\\
V_{40}[-2] = 30381984071/20333356646400 (0.001494)\\
V_{40}[-1] = 2957102734199/20333356646400 (0.145431)\\
V_{40}[+1] = 3360988664941/10166678323200 (0.330589)\\
V_{40}[+2] = 6478070507/10166678323200 (0.000637)\\
V_{40}[+3] = 292144933/6777785548800 (0.000043)\\
V_{40}[+4] = 9938893/5083339161600 (0.000002)\\
V_{40}[+5] = 4597979/20333356646400 (0.000000)\\
V_{40}[+6] = 13/30576476160 (0.000000)\\
\\
e_{41} = 7313311/165473280 (0.044196) [2321237890621440]\\
V_{41}[-6] = 61468511/512400587489280 (0.000000)\\
V_{41}[-5] = 905963641/183000209817600 (0.000005)\\
V_{41}[-4] = 739141561/128100146872320 (0.000006)\\
V_{41}[-3] = 95451411817/427000489574400 (0.000224)\\
V_{41}[-2] = 184063646579/73200083927040 (0.002515)\\
V_{41}[-1] = 3161605893907487/26260530108825600 (0.120394)\\
V_{41}[+1] = 8685630577953149/52521060217651200 (0.165374)\\
V_{41}[+2] = 2865735759811/1281001468723200 (0.002237)\\
V_{41}[+3] = 65007805063/427000489574400 (0.000152)\\
V_{41}[+4] = 12935415793/2562002937446400 (0.000005)\\
V_{41}[+5] = 46277477/18300020981760 (0.000003)\\
V_{41}[+6] = 77476183/1281001468723200 (0.000000)\\
\\
e_{42} = 6123037/27578880 (0.222019) [11945083135426560]\\
V_{42}[-6] = 396751/30500034969600 (0.000000)\\
V_{42}[-5] = 14830471/20333356646400 (0.000001)\\
V_{42}[-4] = 1896873157/142333496524800 (0.000013)\\
V_{42}[-3] = 3201729511/13343765299200 (0.000240)\\
V_{42}[-2] = 33789705539/18977799536640 (0.001780)\\
V_{42}[-1] = 93023983380571/854000979148800 (0.108927)\\
V_{42}[+1] = 283373054708011/854000979148800 (0.331818)\\
V_{42}[+2] = 705994382021/427000489574400 (0.001653)\\
V_{42}[+3] = 6305790613/122000139878400 (0.000052)\\
V_{42}[+4] = 338310517/427000489574400 (0.000001)\\
V_{42}[+5] = 159624961/854000979148800 (0.000000)\\
V_{42}[+6] = 7585/1265186635776 (0.000000)\\
\\
e_{43} = 1241881/47278080 (0.026268) [63663558680494080]\\
V_{43}[-7] = 899/20600900812800 (0.000000)\\
V_{43}[-6] = 239374189/56364064623820800 (0.000000)\\
V_{43}[-5] = 8097232133/56364064623820800 (0.000000)\\
V_{43}[-4] = 35174564923/12525347694182400 (0.000003)\\
V_{43}[-3] = 38747727775/501013907767296 (0.000077)\\
V_{43}[-2] = 34634252980079/18788021541273600 (0.001843)\\
V_{43}[-1] = 779722198548416059/4847309557648588800 (0.160857)\\
V_{43}[+1] = 1906106119849907/10204862226628608 (0.186784)\\
V_{43}[+2] = 8259199022633/4175115898060800 (0.001978)\\
V_{43}[+3] = 536250514949/5368006154649600 (0.000100)\\
V_{43}[+4] = 2687713733/751520861650944 (0.000004)\\
V_{43}[+5] = 6777068293/56364064623820800 (0.000000)\\
V_{43}[+6] = 281711701/112728129247641600 (0.000000)\\
V_{43}[+7] = 899/41201801625600 (0.000000)\\
\\
e_{44} = 27117469/151683840 (0.178776) [443368483462840320]\\
V_{44}[-7] = 19499/1569592442880 (0.000000)\\
V_{44}[-6] = 483086393/521889487257600 (0.000001)\\
V_{44}[-5] = 11256275861/596445128294400 (0.000019)\\
V_{44}[-4] = 664972869923/7515208616509440 (0.000088)\\
V_{44}[-3] = 1013158170881/12525347694182400 (0.000081)\\
V_{44}[-2] = 35958299894189/12525347694182400 (0.002871)\\
V_{44}[-1] = 209668413473641/2087557949030400 (0.100437)\\
V_{44}[+1] = 279518297968507/988843239014400 (0.282672)\\
V_{44}[+2] = 2998032785137/2087557949030400 (0.001436)\\
V_{44}[+3] = 1257623773787/37576043082547200 (0.000033)\\
V_{44}[+4] = 64724291587/37576043082547200 (0.000002)\\
V_{44}[+5] = 50451869/1977686478028800 (0.000000)\\
V_{44}[+6] = 65231/357867076976640 (0.000000)\\
\\
e_{45} = 3979615/45505152 (0.087454) [221817245368320000]\\
V_{45}[-7] = 33263/18788021541273600 (0.000000)\\
V_{45}[-6] = 566438617/10248011749785600 (0.000000)\\
V_{45}[-5] = 269113485517/112728129247641600 (0.000002)\\
V_{45}[-4] = 12278820607/4697005385318400 (0.000003)\\
V_{45}[-3] = 151929636707/1503041723301888 (0.000101)\\
V_{45}[-2] = 46254490783673/28182032311910400 (0.001641)\\
V_{45}[-1] = 6148678815334457/37576043082547200 (0.163633)\\
V_{45}[+1] = 295732896088667/1174251346329600 (0.251848)\\
V_{45}[+2] = 149807819591587/112728129247641600 (0.001329)\\
V_{45}[+3] = 2209747042441/37576043082547200 (0.000059)\\
V_{45}[+4] = 340486901/191714505523200 (0.000002)\\
V_{45}[+5] = 778779179/644160738557952 (0.000001)\\
V_{45}[+6] = 314248783/11272812924764160 (0.000000)\\
V_{45}[+7] = 33263/37576043082547200 (0.000000)\\
\\
e_{46} = 1634903647/27909826560 (0.058578) [151878101868380160]\\
V_{46}[-7] = 593147/1404606873600 (0.000000)\\
V_{46}[-6] = 8041/200658124800 (0.000000)\\
V_{46}[-5] = 461160835213/37576043082547200 (0.000012)\\
V_{46}[-4] = 87404403091/18788021541273600 (0.000005)\\
V_{46}[-3] = 772557276649/18788021541273600 (0.000041)\\
V_{46}[-2] = 14774599478917/6262673847091200 (0.002359)\\
V_{46}[-1] = 30866961312674209/172849798179717120 (0.178577)\\
V_{46}[+1] = 1461674766295439/6173207077847040 (0.236777)\\
V_{46}[+2] = 22857774451559/9394010770636800 (0.002433)\\
V_{46}[+3] = 533933483887/4175115898060800 (0.000128)\\
V_{46}[+4] = 218641508917/18788021541273600 (0.000012)\\
V_{46}[+5] = 2284159373/2087557949030400 (0.000001)\\
V_{46}[+6] = 147831841/2684003077324800 (0.000000)\\
V_{46}[+7] = 8611853/7515208616509440 (0.000000)\\
V_{46}[+8] = 19/2809213747200 (0.000000)\\
\\
e_{47} = 64767749/1213470720 (0.053374) [6786891421192028160]\\
V_{75}[-9] = 1147/715734153953280 (0.000000)\\
V_{75}[-8] = 56233703/169092193871462400 (0.000000)\\
V_{75}[-7] = 8196769949/676368775485849600 (0.000000)\\
V_{75}[-6] = 231016391/27606888795340800 (0.000000)\\
V_{75}[-5] = 115585824097/216438008155471872 (0.000001)\\
V_{75}[-4] = 32560081364423/2705475101943398400 (0.000012)\\
V_{75}[-3] = 31938251933447/193248221567385600 (0.000165)\\
V_{75}[-2] = 10645199114119417/5410950203886796800 (0.001967)\\
V_{75}[-1] = 27034488309751213037/254314659582679449600 (0.106303)\\
V_{75}[+1] = 115130557033516291/722484828359884800 (0.159354)\\
V_{75}[+2] = 5876022667634491/2705475101943398400 (0.002172)\\
V_{75}[+3] = 31091772610243/216438008155471872 (0.000144)\\
V_{75}[+4] = 3905889918467/541095020388679680 (0.000007)\\
V_{75}[+5] = 191742884821/676368775485849600 (0.000000)\\
V_{75}[+6] = 3669340891/772992886269542400 (0.000000)\\
V_{75}[+7] = 32823296417/5410950203886796800 (0.000000)\\
V_{75}[+8] = 56233703/338184387742924800 (0.000000)\\
V_{75}[+9] = 1147/1431468307906560 (0.000000)\\
\\
e_{48} = 3288940867/12134707200 (0.271036) [35197477704016330752]\\
V_{48}[-7] = 25505321/94928950945382400 (0.000000)\\
V_{48}[-6] = 8172767/645775176499200 (0.000000)\\
V_{48}[-5] = 1007637241/47464475472691200 (0.000000)\\
V_{48}[-4] = 194589824999/163968187996569600 (0.000001)\\
V_{48}[-3] = 23433243267407/901825033981132800 (0.000026)\\
V_{48}[-2] = 181000587404083/150304172330188800 (0.001204)\\
V_{48}[-1] = 174759438479524523/1803650067962265600 (0.096892)\\
V_{48}[+1] = 109207543272255949/300608344660377600 (0.363288)\\
V_{48}[+2] = 5722102669224553/1803650067962265600 (0.003173)\\
V_{48}[+3] = 42223142910659/180365006796226560 (0.000234)\\
V_{48}[+4] = 4708918459513/257664295423180800 (0.000018)\\
V_{48}[+5] = 18463193263/9394010770636800 (0.000002)\\
V_{48}[+6] = 23339473/303644792586240 (0.000000)\\
V_{48}[+7] = 1193860483/66801854368972800 (0.000000)\\
V_{48}[+8] = 8671/10844456878080 (0.000000)\\
V_{48}[+9] = 519961/54656062665523200 (0.000000)\\
\\
e_{49} = 10989331/260029440 (0.042262) [5602583726551203840]\\
V_{49}[-7] = 131153887/901825033981132800 (0.000000)\\
V_{49}[-6] = 8877334111/676368775485849600 (0.000000)\\
V_{49}[-5] = 321058535947/676368775485849600 (0.000000)\\
V_{49}[-4] = 188187460081/17568020142489600 (0.000011)\\
V_{49}[-3] = 149860268269/1294485694709760 (0.000116)\\
V_{49}[-2] = 9555417316875511/5410950203886796800 (0.001766)\\
V_{49}[-1] = 77997264836317111/541095020388679680 (0.144147)\\
V_{49}[+1] = 12026849533501/64724284735488 (0.185817)\\
V_{49}[+2] = 5202805469811217/2705475101943398400 (0.001923)\\
V_{49}[+3] = 106923208106387/541095020388679680 (0.000198)\\
V_{49}[+4] = 24986800232473/1352737550971699200 (0.000018)\\
V_{49}[+5] = 2087822266097/2705475101943398400 (0.000001)\\
V_{49}[+6] = 87856669279/5410950203886796800 (0.000000)\\
V_{49}[+7] = 8813377/64416073855795200 (0.000000)\\
\\
e_{50} = 379369/2996224 (0.126616) [17127781233131520000]\\
V_{50}[-7] = 2101783/47464475472691200 (0.000000)\\
V_{50}[-6] = 164943659/33400927184486400 (0.000000)\\
V_{50}[-5] = 465186406151/1803650067962265600 (0.000000)\\
V_{50}[-4] = 49542317623/9109343777587200 (0.000005)\\
V_{50}[-3] = 6825915940891/100202781553459200 (0.000068)\\
V_{50}[-2] = 3811750874427449/1803650067962265600 (0.002113)\\
V_{50}[-1] = 29887357972193783/225456258495283200 (0.132564)\\
V_{50}[+1] = 6548160578576831/25050695388364800 (0.261396)\\
V_{50}[+2] = 708473647746719/901825033981132800 (0.000786)\\
V_{50}[+3] = 356805154396393/1803650067962265600 (0.000198)\\
V_{50}[+4] = 11564370480493/901825033981132800 (0.000013)\\
V_{50}[+5] = 893786131687/225456258495283200 (0.000004)\\
V_{50}[+6] = 468039803071/1803650067962265600 (0.000000)\\
V_{50}[+7] = 104062907/25050695388364800 (0.000000)\\
V_{50}[+8] = 59551/12269728353484800 (0.000000)\\
\\
e_{51} = 12998582579/206290022400 (0.063011) [8694237351200292864]\\
V_{51}[-8] = 302799841/71196713209036800 (0.000000)\\
V_{51}[-7] = 200198636767/1082190040777359360 (0.000000)\\
V_{51}[-6] = 4232845423/28478685283614720 (0.000000)\\
V_{51}[-5] = 296867080061/110427555181363200 (0.000003)\\
V_{51}[-4] = 28279018507781/2705475101943398400 (0.000010)\\
V_{51}[-3] = 5285026941479/48312055391846400 (0.000109)\\
V_{51}[-2] = 2400296876091767/1352737550971699200 (0.001774)\\
V_{51}[-1] = 36173757934527517/238925073937858560 (0.151402)\\
V_{51}[+1] = 9847006606815082297/45993076733037772800 (0.214098)\\
V_{51}[+2] = 9792585533485231/5410950203886796800 (0.001810)\\
V_{51}[+3] = 67216527293609/386496443134771200 (0.000174)\\
V_{51}[+4] = 105773214248153/5410950203886796800 (0.000020)\\
V_{51}[+5] = 1496553336277/284786852836147200 (0.000005)\\
V_{51}[+6] = 1592766194951/5410950203886796800 (0.000000)\\
V_{51}[+7] = 10213095919/27606888795340800 (0.000000)\\
V_{51}[+8] = 302799841/35598356604518400 (0.000000)\\
\\
e_{52} = 4136435441/39437798400 (0.104885) [14755722292396818432]\\
V_{52}[-7] = 154583/1078738078924800 (0.000000)\\
V_{52}[-6] = 8275709/89894839910400 (0.000000)\\
V_{52}[-5] = 5082975040331/601216689320755200 (0.000008)\\
V_{52}[-4] = 4086365658149/200405563106918400 (0.000020)\\
V_{52}[-3] = 146925185217599/450912516990566400 (0.000326)\\
V_{52}[-2] = 3225891061904431/901825033981132800 (0.003577)\\
V_{52}[-1] = 3629547434141844893/23447450883509452800 (0.154795)\\
V_{52}[+1] = 6208463813049931387/23447450883509452800 (0.264782)\\
V_{52}[+2] = 2533165683241/2109532243230720 (0.001201)\\
V_{52}[+3] = 188653240419367/901825033981132800 (0.000209)\\
V_{52}[+4] = 3227651747927/112728129247641600 (0.000029)\\
V_{52}[+5] = 3620914111/1937325529497600 (0.000002)\\
V_{52}[+6] = 233961671281/1803650067962265600 (0.000000)\\
V_{52}[+7] = 3363306073/150304172330188800 (0.000000)\\
V_{52}[+8] = 35069/28757175828480 (0.000000)\\
V_{52}[+9] = 218671/7910745912115200 (0.000000)\\
V_{52}[+10] = 10013/47464475472691200 (0.000000)\\
\\
e_{53} = 1753484317/36404121600 (0.048167) [372961927410354487296]\\
V_{53}[-13] = 1/4831205539184640 (0.000000)\\
V_{53}[-12] = 137923/3043659489686323200 (0.000000)\\
V_{53}[-11] = 70055/15218297448431616 (0.000000)\\
V_{53}[-10] = 24043013971/97397103669962342400 (0.000000)\\
V_{53}[-9] = 536992451/81984093998284800 (0.000000)\\
V_{53}[-8] = 495270318673/7304782775247175680 (0.000000)\\
V_{53}[-7] = 960637729199/24349275917490585600 (0.000000)\\
V_{53}[-6] = 214231483919113/292191311009887027200 (0.000001)\\
V_{53}[-5] = 56397449178799/13281423227722137600 (0.000004)\\
V_{53}[-4] = 3061411267633/216438008155471872 (0.000014)\\
V_{53}[-3] = 15076596535394231/97397103669962342400 (0.000155)\\
V_{53}[-2] = 250559460613448869/146095655504943513600 (0.001715)\\
V_{53}[-1] = 212226881850682202239/1935767435440501555200 (0.109634)\\
V_{53}[+1] = 243813332843085582721/1548613948352401244160 (0.157440)\\
V_{53}[+2] = 105517436829081911/58438262201977405440 (0.001806)\\
V_{53}[+3] = 9720451206859513/48698551834981171200 (0.000200)\\
V_{53}[+4] = 332885526845159/16232850611660390400 (0.000021)\\
V_{53}[+5] = 2199855412954037/292191311009887027200 (0.000008)\\
V_{53}[+6] = 393376392231287/292191311009887027200 (0.000001)\\
V_{53}[+7] = 1765246803409/24349275917490585600 (0.000000)\\
V_{53}[+8] = 10420346053/76892450265759744 (0.000000)\\
V_{53}[+9] = 536992451/40992046999142400 (0.000000)\\
V_{53}[+10] = 24043013971/48698551834981171200 (0.000000)\\
V_{53}[+11] = 70055/7609148724215808 (0.000000)\\
V_{53}[+12] = 137923/1521829744843161600 (0.000000)\\
V_{53}[+13] = 1/2415602769592320 (0.000000)\\
\\
e_{54} = 32744383961/169885900800 (0.192743) [1520584461002805608448]\\
V_{54}[-10] = 519961/2951427383938252800 (0.000000)\\
V_{54}[-9] = 254449/10540812085493760 (0.000000)\\
V_{54}[-8] = 544051261/491904563989708800 (0.000000)\\
V_{54}[-7] = 17186911411/885428215181475840 (0.000000)\\
V_{54}[-6] = 267372853813/1545985772539084800 (0.000000)\\
V_{54}[-5] = 9682910835401/3607300135924531200 (0.000003)\\
V_{54}[-4] = 10292911323827/885428215181475840 (0.000012)\\
V_{54}[-3] = 11916594797204821/97397103669962342400 (0.000122)\\
V_{54}[-2] = 84872899890985/61839430901563392 (0.001372)\\
V_{54}[-1] = 10174542900243355387/97397103669962342400 (0.104465)\\
V_{54}[+1] = 7257687613690421939/24349275917490585600 (0.298066)\\
V_{54}[+2] = 32905496396048207/32465701223320780800 (0.001014)\\
V_{54}[+3] = 294333194742383/4427141075907379200 (0.000066)\\
V_{54}[+4] = 630305608616617/32465701223320780800 (0.000019)\\
V_{54}[+5] = 378893651143/220855110362726400 (0.000002)\\
V_{54}[+6] = 914315912029/2563081675525324800 (0.000000)\\
V_{54}[+7] = 1058275765039/32465701223320780800 (0.000000)\\
V_{54}[+8] = 10007111909/4058212652915097600 (0.000000)\\
V_{54}[+9] = 205655749/2563081675525324800 (0.000000)\\
V_{54}[+10] = 4511477/366154525075046400 (0.000000)\\
V_{54}[+11] = 2800937/5126163351050649600 (0.000000)\\
V_{54}[+12] = 33337/5126163351050649600 (0.000000)\\
\\
e_{55} = 943387771/18202060800 (0.051829) [416456526393846005760]\\
V_{55}[-10] = 33337/8321693751705600 (0.000000)\\
V_{55}[-9] = 83722897/295142738393825280 (0.000000)\\
V_{55}[-8] = 1149868206293/146095655504943513600 (0.000000)\\
V_{55}[-7] = 1532007576709/9130978469058969600 (0.000000)\\
V_{55}[-6] = 50123746848421/97397103669962342400 (0.000001)\\
V_{55}[-5] = 9301249107517/1739233994106470400 (0.000005)\\
V_{55}[-4] = 967852006786867/58438262201977405440 (0.000017)\\
V_{55}[-3] = 6951748079365363/73047827752471756800 (0.000095)\\
V_{55}[-2] = 36857804567275229/29219131100988702720 (0.001261)\\
V_{55}[-1] = 43645307985941263397/292191311009887027200 (0.149372)\\
V_{55}[+1] = 3660220040113447207/18261956938117939200 (0.200429)\\
V_{55}[+2] = 442368372384129727/292191311009887027200 (0.001514)\\
V_{55}[+3] = 489893102204383/2981543989896806400 (0.000164)\\
V_{55}[+4] = 7068099505148281/292191311009887027200 (0.000024)\\
V_{55}[+5] = 10499427847969/1014553163228774400 (0.000010)\\
V_{55}[+6] = 888556569941/869616997053235200 (0.000001)\\
V_{55}[+7] = 98029940871487/292191311009887027200 (0.000000)\\
V_{55}[+8] = 4599455666459/292191311009887027200 (0.000000)\\
V_{55}[+9] = 83722897/147571369196912640 (0.000000)\\
V_{55}[+10] = 33337/4160846875852800 (0.000000)\\
\\
e_{56} = 400879067/2316625920 (0.173044) [1415737696661966684160]\\
V_{56}[-6] = 6785814997/10821900407773593600 (0.000000)\\
V_{56}[-5] = 10357492793639/97397103669962342400 (0.000000)\\
V_{56}[-4] = 161442370258651/97397103669962342400 (0.000002)\\
V_{56}[-3] = 181501819795279/6493140244664156160 (0.000028)\\
V_{56}[-2] = 87969285921507127/97397103669962342400 (0.000903)\\
V_{56}[-1] = 66389099956276679/512616335105064960 (0.129510)\\
V_{56}[+1] = 2903737047336093503/9739710366996234240 (0.298134)\\
V_{56}[+2] = 13339984584863687/5410950203886796800 (0.002465)\\
V_{56}[+3] = 9234838162290949/24349275917490585600 (0.000379)\\
V_{56}[+4] = 88940575640699/1739233994106470400 (0.000051)\\
V_{56}[+5] = 10314251483003/1475713691969126400 (0.000007)\\
V_{56}[+6] = 75172645647961/48698551834981171200 (0.000002)\\
V_{56}[+7] = 9251203375/68348844680675328 (0.000000)\\
V_{56}[+8] = 15404457773/2705475101943398400 (0.000000)\\
V_{56}[+9] = 4004766797/32465701223320780800 (0.000000)\\
V_{56}[+10] = 903919/507276581614387200 (0.000000)\\
V_{56}[+11] = 403/26698767453388800 (0.000000)\\
V_{56}[+12] = 41/2029106326457548800 (0.000000)\\
\\
e_{57} = 30907469237/509657702400 (0.060644) [505006549580537266176]\\
V_{57}[-9] = 9899827187/292191311009887027200 (0.000000)\\
V_{57}[-8] = 2078227018577/292191311009887027200 (0.000000)\\
V_{57}[-7] = 279694033471/4869855183498117120 (0.000000)\\
V_{57}[-6] = 57975124388297/146095655504943513600 (0.000000)\\
V_{57}[-5] = 531949181968363/146095655504943513600 (0.000004)\\
V_{57}[-4] = 31272327932357/1987695993264537600 (0.000016)\\
V_{57}[-3] = 12585903870578537/97397103669962342400 (0.000129)\\
V_{57}[-2] = 57593752652897383/32465701223320780800 (0.001774)\\
V_{57}[-1] = 43587023408027517991/292191311009887027200 (0.149173)\\
V_{57}[+1] = 1530357192955689677/7304782775247175680 (0.209501)\\
V_{57}[+2] = 44651341400352107/24349275917490585600 (0.001834)\\
V_{57}[+3] = 777758638132831/4058212652915097600 (0.000192)\\
V_{57}[+4] = 148197040200331/10821900407773593600 (0.000014)\\
V_{57}[+5] = 1954397988395299/292191311009887027200 (0.000007)\\
V_{57}[+6] = 182673651588179/292191311009887027200 (0.000001)\\
V_{57}[+7] = 1351993586243/12174637958745292800 (0.000000)\\
V_{57}[+8] = 2078227018577/146095655504943513600 (0.000000)\\
V_{57}[+9] = 9899827187/146095655504943513600 (0.000000)\\
\\
e_{58} = 6267235837/70718054400 (0.088623) [750950012179833421824]\\
V_{58}[-9] = 324250979/2563081675525324800 (0.000000)\\
V_{58}[-8] = 23369547097/2563081675525324800 (0.000000)\\
V_{58}[-7] = 83894965921/512616335105064960 (0.000000)\\
V_{58}[-6] = 287011511657/732309050150092800 (0.000000)\\
V_{58}[-5] = 767530296899437/97397103669962342400 (0.000008)\\
V_{58}[-4] = 74756283949183/8854282151814758400 (0.000008)\\
V_{58}[-3] = 9067008775798753/97397103669962342400 (0.000093)\\
V_{58}[-2] = 122681701888511237/97397103669962342400 (0.001260)\\
V_{58}[-1] = 189447703516884708733/1412258003214453964800 (0.134145)\\
V_{58}[+1] = 207458630795906734361/941505335476302643200 (0.220348)\\
V_{58}[+2] = 18651227966840171/8116425305830195200 (0.002298)\\
V_{58}[+3] = 2862258848807693/16232850611660390400 (0.000176)\\
V_{58}[+4] = 573932994775373/16232850611660390400 (0.000035)\\
V_{58}[+5] = 11734282833971/2434927591749058560 (0.000005)\\
V_{58}[+6] = 26878790210113/48698551834981171200 (0.000001)\\
V_{58}[+7] = 23257660896641/97397103669962342400 (0.000000)\\
V_{58}[+8] = 215069764351/19479420733992468480 (0.000000)\\
V_{58}[+9] = 230825334239/48698551834981171200 (0.000000)\\
V_{58}[+10] = 5755065131/24349275917490585600 (0.000000)\\
V_{58}[+11] = 68623789/24349275917490585600 (0.000000)\\$
\\
\vfil\eject
\noindent $e_{59} = 163742997419/3057946214400 (0.053547) [27693220083749128765440]\\
V_{59}[-15] = 41/6763687754858496000 (0.000000)\\
V_{59}[-14] = 61/84546096935731200 (0.000000)\\
V_{59}[-13] = 1319/48312055391846400 (0.000000)\\
V_{59}[-12] = 2483197/6763687754858496000 (0.000000)\\
V_{59}[-11] = 4760227907/1460956555049435136000 (0.000000)\\
V_{59}[-10] = 36261195821/292191311009887027200 (0.000000)\\
V_{59}[-9] = 99251424791/91309784690589696000 (0.000000)\\
V_{59}[-8] = 145306469550229/5843826220197740544000 (0.000000)\\
V_{59}[-7] = 96842410787587/922709403189116928000 (0.000000)\\
V_{59}[-6] = 908924421985849/834832317171105792000 (0.000001)\\
V_{59}[-5] = 15866293406797/2385235191917445120 (0.000007)\\
V_{59}[-4] = 7538369215801399/273929354071769088000 (0.000028)\\
V_{59}[-3] = 469300376053615651/1593770787326656512000 (0.000294)\\
V_{59}[-2] = 30998833972906813663/17531478660593221632000 (0.001768)\\
V_{59}[-1] = 101336333823239581902691/1034357240975000076288000 (0.097970)\\
V_{59}[+1] = 15731441434639484641049/103435724097500007628800 (0.152089)\\
V_{59}[+2] = 503906908616282993/318754157465331302400 (0.001581)\\
V_{59}[+3] = 2029679255343050999/8765739330296610816000 (0.000232)\\
V_{59}[+4] = 17982332246708579/626124237878329344000 (0.000029)\\
V_{59}[+5] = 453773605688981/106251385821777100800 (0.000004)\\
V_{59}[+6] = 147689548546613/166966463434221158400 (0.000001)\\
V_{59}[+7] = 1221525888512777/17531478660593221632000 (0.000000)\\
V_{59}[+8] = 116187185196623/5843826220197740544000 (0.000000)\\
V_{59}[+9] = 5560595677129/5843826220197740544000 (0.000000)\\
V_{59}[+10] = 147844102799/1460956555049435136000 (0.000000)\\
V_{59}[+11] = 1769351693/417416158585552896000 (0.000000)\\
V_{59}[+12] = 2483197/3381843877429248000 (0.000000)\\
V_{59}[+13] = 1319/24156027695923200 (0.000000)\\
V_{59}[+14] = 61/42273048467865600 (0.000000)\\
V_{59}[+15] = 41/3381843877429248000 (0.000000)\\$
\\
\vfil\eject
\noindent $e_{60} = 67410047209/254828851200 (0.264531) [139128414846556255027200]\\
V_{60}[-11] = 48572009/307569801063038976000 (0.000000)\\
V_{60}[-10] = 7564923041/102523267021012992000 (0.000000)\\
V_{60}[-9] = 241946098633/61513960212607795200 (0.000000)\\
V_{60}[-8] = 82710948197627/5843826220197740544000 (0.000000)\\
V_{60}[-7] = 618201928419193/1460956555049435136000 (0.000000)\\
V_{60}[-6] = 3617639494318727/1460956555049435136000 (0.000002)\\
V_{60}[-5] = 25317575928648679/1947942073399246848000 (0.000013)\\
V_{60}[-4] = 56398540189930703/531256929108885504000 (0.000106)\\
V_{60}[-3] = 93326831348910953/584382622019774054400 (0.000160)\\
V_{60}[-2] = 2254270194640917743/834832317171105792000 (0.002700)\\
V_{60}[-1] = 63503355003238436449/649314024466415616000 (0.097801)\\
V_{60}[+1] = 2119164534676506715771/5843826220197740544000 (0.362633)\\
V_{60}[+2] = 2258443180589756309/834832317171105792000 (0.002705)\\
V_{60}[+3] = 105935771427821093/531256929108885504000 (0.000199)\\
V_{60}[+4] = 540048225702181/39753919865290752000 (0.000014)\\
V_{60}[+5] = 3885941279688631/973971036699623424000 (0.000004)\\
V_{60}[+6] = 7952143409509/26088509911597056000 (0.000000)\\
V_{60}[+7] = 68408181635527/973971036699623424000 (0.000000)\\
V_{60}[+8] = 119925986148293/5843826220197740544000 (0.000000)\\
V_{60}[+9] = 5490125991767/5843826220197740544000 (0.000000)\\
V_{60}[+10] = 244847740997/307569801063038976000 (0.000000)\\
V_{60}[+11] = 14716614571/973971036699623424000 (0.000000)\\
V_{60}[+12] = 2237256377/243492759174905856000 (0.000000)\\
V_{60}[+13] = 30677699/177085643036295168000 (0.000000)\\$
\\
\vfil\eject
\noindent $e_{61} = 19902302711/764486553600 (0.026034) [863065094810581724037120]\\
V_{61}[-13] = 33337/19223112566439936000 (0.000000)\\
V_{61}[-12] = 26683982221/18115861282612995686400 (0.000000)\\
V_{61}[-11] = 4775116095173/30193102137688326144000 (0.000000)\\
V_{61}[-10] = 69756035473169/15096551068844163072000 (0.000000)\\
V_{61}[-9] = 312279556434773/45289653206532489216000 (0.000000)\\
V_{61}[-8] = 72308884322139191/362317225652259913728000 (0.000000)\\
V_{61}[-7] = 34398025829829313/120772408550753304576000 (0.000000)\\
V_{61}[-6] = 457539002368786193/543475838478389870592000 (0.000001)\\
V_{61}[-5] = 2575014867714166963/181158612826129956864000 (0.000014)\\
V_{61}[-4] = 3066031235443749497/155278810993825677312000 (0.000020)\\
V_{61}[-3] = 15639813937655487379/67934479809798733824000 (0.000230)\\
V_{61}[-2] = 1423557203480883778099/1086951676956779741184000 (0.001310)\\
V_{61}[-1] = 2061939117042950979905/13958747851444960886784 (0.147717)\\
V_{61}[+1] = 11468365862330930545275961/66304052294363564212224000 (0.172966)\\
V_{61}[+2] = 44062736926197420397/24703447203563175936000 (0.001784)\\
V_{61}[+3] = 25978409273977745917/135868959619597467648000 (0.000191)\\
V_{61}[+4] = 19021961155258063699/1086951676956779741184000 (0.000018)\\
V_{61}[+5] = 89329015884779113/12077240855075330457600 (0.000007)\\
V_{61}[+6] = 2382219017953771/5200725727065931776000 (0.000000)\\
V_{61}[+7] = 4359398737953311/25879801832304279552000 (0.000000)\\
V_{61}[+8] = 185240586402109/1848557273736019968000 (0.000000)\\
V_{61}[+9] = 312280878311483/90579306413064978432000 (0.000000)\\
V_{61}[+10] = 69756035473169/30193102137688326144000 (0.000000)\\
V_{61}[+11] = 4775116095173/60386204275376652288000 (0.000000)\\
V_{61}[+12] = 26683982221/36231722565225991372800 (0.000000)\\
V_{61}[+13] = 33337/38446225132879872000 (0.000000)\\$
\\

\vfil\eject

\medskip\noindent {\footnotesize Bestevaerstraat 46,
1056 HP Amsterdam, The Netherlands.\\
e-mail: {\tt hjhom48@hotmail.com}}

\medskip\noindent {\footnotesize KdV Institute,
Plantage Muidergracht 24, 1018 TV Amsterdam, The Netherlands.\\
e-mail: {\tt moree@science.uva.nl} (to whom correspondence should be addressed)}

\end{document}